\documentclass[11pt, amssymb,amsfonts]{amsart}
\usepackage{amsrefs}
\usepackage{algorithm}
\usepackage{algorithmic}
\usepackage{subcaption} 

\makeatletter
\def\l@section{\@tocline{1}{10pt plus0pt}{0pt}{}{\bfseries}}

\def\@tocline#1#2#3#4#5#6#7{\relax
    \ifnum #1>-1
  \ifnum #1>\c@tocdepth 
  \else
    \par \addpenalty\@secpenalty
    \begingroup \hyphenpenalty\@M
    \@ifempty{#4}{%
      \@tempdima\csname r@tocindent\number#1\endcsname\relax
    }{%
      \@tempdima#4\relax
    }%
    \parindent\z@ \leftskip#3\relax \advance\leftskip\@tempdima\relax
    \rightskip\@pnumwidth plus4em \parfillskip-\@pnumwidth
    #5\leavevmode\hskip-\@tempdima #6\nobreak\relax
    \hfil\hbox to\@pnumwidth{\@tocpagenum{#7}}\par
    \nobreak
    \endgroup
  \fi
\fi}

\makeatother

\usepackage{a4wide}
\usepackage[a4paper, margin=2.3cm]{geometry}
\usepackage{amsmath,amssymb,amsthm}
\usepackage{color}
\usepackage{parskip}
\usepackage{graphicx}
\usepackage{hyperref}
\usepackage{parskip}
\usepackage{setspace}
\usepackage{amsmath}
\usepackage{mathtools}
\usepackage{amssymb}
\usepackage{cases}
\usepackage{algorithm}
\usepackage{algorithmic} 
\usepackage{todonotes}
\usepackage{stmaryrd}
\usepackage{IEEEtrantools}
\usepackage[nocompress]{cite}
\usepackage{tikz}
\usetikzlibrary{patterns,snakes}
\newcommand{\norm}[1]{ \left|  #1 \right| }
\newcommand{\Norm}[1]{ \left\|  #1 \right\| }

\def\P{{\hbox{\bf P}}}
\def\E{{\hbox{\bf E}}}

 at 10 true pt

\def\be#1{ \begin{equation}\label{#1} }

\def\bas{\begin{equation*}}
\def\eas{\end{equation*}}
\def\bi{\begin{itemize}}
\def\ei{\end{itemize}}

\def\emph#1{{\it #1}}
\def\textbf#1{{\bf #1}}

\theoremstyle{plain}
 \theoremstyle{plain}
  \newtheorem{theorem}{Theorem}
  \numberwithin{theorem}{section}
  
  \newtheorem{problem}{Problem}
  \numberwithin{problem}{section}

  \newtheorem{lemma}{Lemma}
  
  \numberwithin{notation}{section}
  \numberwithin{lemma}{section}
  \newtheorem{corollary}{Corollary}
   \numberwithin{corollary}{section}

\theoremstyle{remark}
  \newtheorem{remark}{Remark}
  \numberwithin{remark}{section}

\theoremstyle{definition}
  \newtheorem{definition}{Definition}
\numberwithin{definition}{section}

\theoremstyle{definition}
  
  \numberwithin{summary}{section}

\DeclareUnicodeCharacter{2113}{\ensuremath{\ell}}
\DeclareUnicodeCharacter{2082}{\ensuremath{\infty}}
\DeclareUnicodeCharacter{221E}{\ensuremath{\infty}}
\begin{document}
\include{psfig}
\title{Matrices perturbed by random noise: The accuracy of low-rank approximation}
\pagenumbering{arabic}

\author{Phuc Tran, Van Vu }
\thanks{trandangphuc234@gmail.com,   School of Engineering, Vin University (VinUni) \\ vuhavan001@gmail.com,  Department of Mathematics,  The University of Hong Kong (HKU) } 
\date{}

\begin{abstract} 
 Let $A$ be an $m \times n$ matrix with rank $r$ and singular value decomposition $A = \sum _{i=1}^r \sigma_i u_i v_i^\top, $ where the $\sigma_i$ are its singular values, ordered decreasingly, and $u_i, v_i$ are the corresponding left and right singular vectors. 
For an integer $1 \le p \le r$, $A_p := \sum_{i=1}^p \sigma_i  u_i v_i^\top$ is the best rank-$p$ approximation of $A$.  In practice, one often chooses $p$ to be small, leading to the 
commonly used phrase ``low-rank approximation''.

\vskip2mm

Low-rank approximation plays a central role in data science because it can substantially reduce the dimensionality of the original data. For a large data matrix $A$, one typically computes a rank-$p$ approximation $A_p$  for a suitably chosen small $p$, stores $A_p$, and uses it as input for further computations. The reduced dimension of $A_p$ enables faster computations and significant data compression.

\vskip2mm 

In practice, noise is inevitable. We 
often have access only to noisy data $\tilde A = A + E$,  where $E$ represents the noise. Consequently, the low-rank approximation used as input in many downstream tasks is $\tilde A_p$, the best rank-$p$ approximation of $\tilde A$, rather than $A_p$. Therefore, it is natural and important to estimate the error $ \| \tilde A_p - A_p \|$. This error plays a critical role in assessing the accuracy of downstream procedures involving low-rank approximations of noisy inputs.
 The standard way to estimate this error is to use the Eckart--Young--Mirsky identity, which provides 
 $$\| \tilde A_p -A_p\| = O(\sigma_{p+1} + \|E\|).$$

A situation that occurs frequently in applications is that noise is random and $A$ has relatively low rank. In  this situation, $\sigma_{p+1}$ often dominates $\|E \|$ by a large factor.
Our main results in this paper show that in this situation, it is possible to improve the Eckart--Young--Mirsky bound by removing the term $\sigma_{p+1}$  completely or replacing it by a much smaller quantity. In certain basic settings, our bounds are sharp, as matching lower bounds 
can be obtained from results in random matrix theory. 

\vskip2mm 
Our proof relies on a contour expansion argument and yields a deterministic result. 
This contour argument is new and could be of independent interest. 

\vskip2mm 
These new results have significant potential in applications. They can be applied to improve the analysis of a broad class of algorithmic tasks involving low-rank approximation of noisy inputs.

\vskip3mm

\textbf{Mathematics Subject Classifications: } 47A55, 68W40.

\end{abstract} 
\maketitle

\section{Introduction} \label{sec: intro}

Let $A$ be an $m \times n$ matrix  with rank $r$ and singular value decomposition (SVD) 
$A = \sum_{i=1}^r \sigma_i u_i v_i^\top,$ where $\sigma_i$ are its singular values, ordered decreasingly, and $u_i, v_i$ are the corresponding left and right singular vectors. We denote by $\delta_k$ the singular value gap 
$\sigma_k- \sigma_{k+1}$ for $1\le k \le r$. 
For a parameter $1 \le p \le r$, the  best rank-$p$ approximation of $A$ is 
$$A_p := \sum_{i=1}^p \sigma_i  u_i v_i^\top . $$

Most of the time, we choose $p$ to be small, focusing on a few leading components in the sum  (which leads to the phrase ``low-rank approximation"). 
Computing a low-rank approximation 
is a task of fundamental interest, with applications in many areas of computer science and statistics; see, for instance, \cite{CT1, dwork2006calibrating, MangoubiVJACM, AchlioptasMcSherry2007, Cha1, chaudhuri2013near, upadhyay2018price, KT1, gonem2018smooth, CW1, TranlinhVu, kolda2009tensor, mahoney2009cur, trefethen2022numerical, golub2013matrix, halko2011finding}. 

In many applications, computing a low-rank approximation serves as an intermediate step rather than the end goal.
Given a large data matrix $A$, one forms a rank $p$  approximation ($A_p)$ for a suitably small $p$, stores 
$A_p$, and then uses it as input for subsequent computations. A common practice is to choose the smallest $p$ so that $A_p$ captures most of (say $90\%$) the ``energy" 
 of $A$, namely, 
 $\|A_p \|_F^2= \sum_{i=1}^p \sigma_i^2 \ge 0.9  \| A \| _F^2  =  \sum_{i=1}^r \sigma_i^2; $
%
see \cite{holmes2023introduction, jolliffe2016principal}.  (The chosen percentage may vary from application to application.) The reduced dimension of $A_p$ enables faster processing and substantial data compression.
Another popular use of low-rank approximation is PCA (Principal Component Analysis). Here, researchers often set 
$p=1,2$ or $3$ to visualize the data; see, for instance \cite{jolliffe2016principal, golyandina2014basic, candes2011robust}.

In practice, noise is inevitable. Instead of $A$, one often has access to a noisy matrix  $\tilde A = A+E$, where $E$ represents the noise. The low-rank approximation is then $\tilde A_p$, the best rank-$p$ approximation of $\tilde A$. Thus, in downstream tasks, one works with $\tilde A_p$ as input. In order to evaluate the final output of the whole process, it is natural and important to estimate the error $\| \tilde A_p - A_p \| $. 


{\bf \noindent Main Question. }
{\it Bound $\|\tilde A_p - A_p \|$. } 
\begin{remark} We use the spectral norm $ \| .\| $. Since $\tilde A_p -  A_p$ has rank at most $2p$, we obtain 
a bound in the Frobenius norm by paying an extra factor $\sqrt{2p}$.   
\end{remark}

The standard approach to bounding $\|\tilde A_p-A_p\|$ is to combine the Eckart--Young--Mirsky identity with the triangle inequality. Recall that the Eckart--Young--Mirsky identity states
\(
\|A-A_p\|=\sigma_{p+1};
\)
see \cite{EY1,golub2013matrix}.
\begin{theorem}\label{EY0}
One has
\begin{equation}\label{EKbound}
\|\tilde A_p-A_p\|
\le
2\bigl(\sigma_{p+1}+\|E\|\bigr).
\end{equation}
\end{theorem}
\noindent\textit{Proof.}
By the triangle inequality,
\begin{equation*}
\begin{split}
\|\tilde A_p-A_p\|
&\le
\|\tilde A_p-\tilde A\|
+\|\tilde A-A\|
+\|A-A_p\| \\
&=
\tilde\sigma_{p+1}
+\|E\|
+\sigma_{p+1}
\qquad
\text{(by the Eckart--Young--Mirsky identity)}\\
&\le
2\bigl(\sigma_{p+1}+\|E\|\bigr),
\qquad
\text{(by Weyl's inequality \cite{Book1,We1}).}
\end{split}
\end{equation*}
{\it \noindent Notation.} We use the asymptotic notations 
 $o, O, \Theta,\, etc.,$ under the assumption that $\min \{m, n \} \rightarrow \infty$. For two quantities $ f$ and $g$, we write $f =\tilde O (g)$ if there is a positive constant $c$ such that $f =O( g \log^c \max \{ m,n \}  )$.

 In this paper, we consider a special setting that arises frequently in practice. First, it is very common that noise has a random nature. Next, in many applications, the rank of the ground-truth matrix $A$ is relatively small compared to its dimensions.  We provide four examples of this phenomenon in Appendix \ref{App: motivating examples}, all of which are well-known problems in statistics and theoretical computer science. In these examples, the rank is often  $O(1)$ or $\tilde O(1)$. The small rank phenomenon has also been observed in numerous real-world datasets, and a theoretical explanation for this phenomenon was recently proposed in \cite{UT1}. These two assumptions have been considered in many previous papers which investigate the possibility of improving classical perturbation bounds; see, for instance, \cite{OVK13, TranVishnoiVu2025, ZZ1, CP1, OVK22, EBW1} and the references therein. 

In this paper, we show that in this setting (a relatively low-rank matrix perturbed by random noise), we can obtain a notable improvement upon Theorem \ref{EY0}.  To start, let us observe that between the two terms $\|E \|$ and $\sigma_{p+1}$, the term $\sigma_{p+1}$ is often the dominant one; see the discussion following \eqref{boundonE}. Our study shows that in many situations, we can significantly reduce this term, or even remove it completely. On the other hand, the term $\| E \|$ cannot be improved, due to results from random matrix theory; see Section \ref{sharpness}. 
 
We conclude this introductory discussion by presenting a few corollaries of our main results, which will be presented  in Section \ref{sec: randomnoise}. 

\subsection{A few illustrative results}
\label{newsimplebound}
 Let $A$ be an $m \times n$ matrix of rank $r$. 
Without loss of generality, we assume that $n \ge m$ and  the entries of $A$ have root-mean-square magnitude of order $\Omega (1)$, namely
\begin{equation} \label{Fsum} \sum_{i=1}^r \sigma_i^2 = \| A\|_F^2= \Omega (mn). \end{equation}
 Next, we  let  $E$ be a random matrix whose entries $\xi_{ij}$ 
are independent sub-Gaussian variables with mean 0 and variance 1. A well-known result from random matrix theory shows that with probability $1-o(1) $  \cite{Ver1book},
\begin{equation} \label{boundonE}  \| E \| = \Theta (\sqrt {m+n} ) =\Theta (\sqrt n). \end{equation}

Let $\sigma_{\rm ave}$ be the (square-) average of the singular values, namely $ \sigma_{\rm ave}^2 =\frac{\sum_{i=1}^r \sigma_i^2} {r} $.
By \eqref{Fsum},   $\sigma_{\rm ave} = \Omega (\sqrt{r^{-1} mn}) $. 
As $r \le \min \{m,n\} =m $, it is clear that $\sigma_{\rm ave} $ is at least of the same order of magnitude as 
$\| E \|$. In the setting considered in this paper, where $r$ is  relatively small (say, $r =o(m)$), we would have $
\frac{\sigma_{\rm ave}}{\| E\|} \rightarrow \infty$. 
Furthermore,  in many applications $p=o(1)$ and so $\sigma_{p+1}$ is one of the largest singular values, which is often much bigger than the average $\sigma_{\rm ave}$. In this situation, the term $\sigma_{p+1}$ in \eqref{EKbound} often dominates the term $\| E \|$ by a large factor. If $r$ itself is of order $O(1)$, then 
$\frac{\sigma_{\rm ave}}{\| E \| } = \Omega (\sqrt m)$. 

The following results, which are corollaries of more general theorems (presented later), show that under mild assumptions, we can remove the term  $\sigma_{p+1}$ completely or replace it by a much smaller quantity.  We set $\delta_p := \sigma_p -\sigma_{p+1}$.

\begin{theorem} \label{toy1} Consider the above setting. There is a constant $C > 0$ so that if $\delta_p \ge Cr \log n \cdot \max \big\{ r^2 , \frac{  \sigma_p }{\sqrt{n}}, \frac{n}{\sigma_p} \big\}$,  then with probability $1-o(1)$, 
$$\| \tilde A_p - A_p \| = O( \| E \| ) = O(\sqrt n ). $$
\end{theorem}

Notice that if $\sigma_1 \leq C' \sigma_p$ for some constant $C'$, then $\sigma_p = \Omega( \| A \| _F/\sqrt{r} )  = \Theta(\sqrt {r^{-1}mn} )$. Naturally, we also have $\sigma_p \leq \|A\|_F=O(\sqrt{mn})$.
Thus, we have the following corollary. 

\begin{corollary} \label{toy1-1} Consider the above setting. There is a constant $C > 0$ so that if $\delta_p \ge C r \log n \cdot \max  \big\{ r^2, \sqrt{m},    \sqrt{\frac{rn}{m}}  \big\}$, then with probability $1-o(1)$, 
$$\| \tilde A_p - A_p \| = O( \| E \| ) = O(\sqrt n ). $$
\end{corollary}

In particular, if $r =O(1)$, we can simplify the factor with $r$ and obtain 
\begin{corollary} \label{toy1-2} Consider the above setting. Assume that $r=O(1)$. There is a constant $C > 0$ so that if $\delta_p \ge C  \log n \cdot \max  \big\{ \sqrt{m},    \sqrt{\frac{n}{m}}  \big\}$, then with probability $1-o(1)$, 
$$\| \tilde A_p - A_p \| = O( \| E \| ) = O(\sqrt n ). $$
\end{corollary}

\begin{theorem} \label{toy3} 
Consider the above setting. Assume that 
$r \leq \log^{-1} n \cdot \min\{\sqrt{\delta_p}, \frac{\delta_p \sigma_p}{n}  \} $ and $\sigma_p =\Omega(\sigma_{\mathrm{ave}})$. 
Set $\epsilon_0 := \max \{ \frac{r \log n}{\sqrt \delta_p}, \frac{r^{7/2} \log n} {\delta_p \sqrt m } \} $. Then  with probability $1-o(1)$ 
$$\| \tilde A_p - A_p \| = O( \| E \|  + \epsilon_0 \sigma_{p+1} ). $$

\end{theorem}

 \begin{remark} \label{smallgap}  (Small gap conditions.) Notice that all of these results involve a lower bound on the gap $\delta_p$. Such a requirement is necessary; see Subsection \ref{gap-nec} for more details. 
 Notice that if  $r$ is small, we can have $\delta_p$ much smaller than $\sigma_p$. This means that  $\sigma_p$ and $\sigma_{p+1}$ are allowed to be relatively close to each other. See  Remark \ref{remark: daviskahan}  for a comparison to the gap conditions needed for using the Davis-Kahan theorem. 
\end{remark}

\subsection{Structure of the rest of the paper}  \label{restofpaper}

\begin{itemize}

\item In Section \ref{sec: randomnoise}, we present general results with respect to random noise, which are far-reaching generalizations of the illustrative theorems presented above.  

\item In Section \ref{subsec: mainDe}, we present
 Theorem \ref{cor: rec}, which applies to a deterministic matrix  $E$. We then comment on our new key parameters,  the skewness parameters $x$ and $y$,  and work out a representative example to illustrate their power in our study. Also in this section, we discuss the low-rank condition on $A$. We point out that, in a subsequent paper \cite{DKTranVu2}, we can replace this by 
a weaker assumption that the singular values do not cluster near $\sigma_p$; see Subsection \ref{subsec: notlowrank}.
 We conclude this section with a geometric perspective on the main results.  

\item In Section \ref{section: symmetric}, we present the symmetric version of Theorem \ref{cor: rec}. We can reduce the non-symmetric case to the symmetric one by a simple symmetrization trick.


\item  In Section \ref{sec: proofidea}, we present a sketch of the proof. In the next three sections,  Sections \ref{sec: proof}, \ref{proof: lemma}, and \ref{proof:lemma2}, we give the detailed proof of the main theorem and its supporting lemmas.

\item In Section \ref{sec: Apps}, we discuss an application, in which we use our new results to improve the error analysis of computing low-rank approximation for a matrix with missing and noisy entries.

\end{itemize}
%
%
%
%

\section{ General results with random noise} \label{sec: randomnoise}

In this section, we present general results for random noise. Let $E$ be an $m \times n$ real random matrix whose entries $\xi_{ij}$ are independent (but not necessarily iid) random variables with mean zero and variance  $\sigma_{ij}^2$. 
We say that $E$ is
\begin{itemize} 

\item $K$-bounded if with probability $1$,  $|\xi_{ij}| \le K$, for all pairs $i,j$,

\item $(K,\sigma)$-bounded if with probability $1$, $ |\xi_{ij}| \leq K$ and $ \sigma_{ij}^2 \leq \sigma^2$ for all  pairs $i,j$. 

\end{itemize}

Denote 
$$m_4 := \max_{ 1 \leq k \leq n, 1 \leq l \leq m} \bigg\{ \frac{1}{m }\sum_{i=1}^{m} \sigma_{ik}^4 , \frac{1}{n}\sum_{j=1}^{n} \sigma_{lj}^4 \bigg\}.$$
Thus, $m_4$ is the maximum, over rows and columns, of the average squared variances. 

The following is our key assumption:
\begin{equation} \label{ass: mainran}
 \textbf{(C0):} \,\,\,\,  \max \left\lbrace  p\epsilon_1, r^2 t_1 \epsilon_2 , \sqrt{r \eta} \right\rbrace \leq \frac{1}{96},
\end{equation}
where
\begin{equation} \label{epsilon12-App} \epsilon_1:=\frac{\|E\|}{\sigma_p}   , \epsilon_2:=\frac{1} {\delta_p }, \eta:=\frac{\|E\|^2}{\sigma_p \delta_p }.\end{equation}
Here, the parameter $t_1$ is part of our theorems. The constant $96$ is chosen for convenience and can be considerably reduced by tightening the analysis in  Section~\ref{sec: proof}.

\begin{remark} \label{remark:C0} (A discussion on Assumption {\bf C0}.) As we will use {\bf C0} and similar assumptions in all of our theorems, let us comment on its nature. As $p$ is small, the first condition $p \epsilon_1 \le \frac{1}{96}$
simply means that $\| E\|$ is smaller than $\sigma_p$, which, as discussed above, often holds automatically.  The second condition 
$r^2 t_1 \epsilon_2 \le \frac{1}{96}$ is basically a bound  on the rank $r$, providing that $r$ is small compared to $\sqrt {\delta_p } $. (The parameter $t_1$ is typically of logarithmic order.) Finally, the condition 
$\sqrt {r \eta} \le \frac{1}{96}$ is guaranteed if
$$\delta_p \gg r \frac{\|E \|^2}{ \sigma_p}.$$

 Intuitively, it is expected that some lower bound on $\delta_p$ is 
necessary, as if this is too small, then the noise may swap the two singular values
$\sigma_p$ and $\sigma_{p+1}$ and their corresponding singular vectors, leading to a significant perturbation of $A_p$. In fact, studies in random matrix  theory showed that the ratio 
$\frac{\| E\|^2}{\sigma_p}$ is the correct order of the perturbation of $\sigma_p$. 
Thus, having $\delta_p \gg \frac{\| E \|^2}{\sigma_p}$ would help us to avoid the potential singular value swap. 
We will discuss this point in more detail in Subsection \ref{gap-nec}.  \end{remark}

\begin{theorem} \label{theo: mainTh-1}
Let $A$ be an $m \times n$ matrix of rank $r$  ($m \leq n$) and $E$ be a random matrix, 
whose entries are iid $K$-bounded random variables with zero mean and variance one. 
Assume that  \textbf{C0} holds for some parameter $t_1$.  Then for any $t_2 >0$, we have 
$$\Norm{\tilde{A}_p - A_p} \leq 32 \sigma_p \left(  \epsilon_1 + rt_1 \epsilon_2 + \frac{r^2 t_2 }{\|E\|} \epsilon_1 \epsilon_2 \right) = 32 \left( \|E\|+ r t_1 \frac{\sigma_p}{\delta_p} + \frac{r^2 t_2 }{\delta_p} \right),$$
with probability at least $1 - r^2 \exp\left( \frac{-t_1^2/2}{4 + Kt_1} \right) - \frac{5r(r-1)(m+n) m_4}{t_2^2}. $ 
\end{theorem}

Although Theorem \ref{theo: mainTh-1} appears technical at first, it is straightforward to apply. One first chooses the parameters $t_1$ and $t_2$ to guarantee a desired probability bound. When the rank $r$ is not too large,  both $t_1$ and $\frac{t_2}{\sqrt{n}}$ will be small, making the error terms $ r t_1 \frac{\sigma_p}{\delta_p} $ and $ \frac{r^2 t_2 }{\delta_p}$  both  small. 


We illustrate this strategy by proving Theorem~\ref{toy1}. Here, $m_4=1$, and a standard truncation argument (see the next subsection) allows us to take $K=O(\log n)$. We can set 
\[
t_1=10K\log n,
\qquad
t_2=r\sqrt n\log n,
\]
to guarantee that the probability bound is $1-o(1)$. By the assumption on $\delta_p$ in Theorem \ref{toy1}, we have 
\[ \textstyle
rt_1\frac{\sigma_p}{\delta_p} = O\big(\frac{r\sigma_p\log n}{\delta_p}\big) = O(\sqrt n) =
O(\|E\|),
\]
and similarly,
\[ \textstyle
\frac{r^2t_2}{\delta_p} = O\big(\frac{r^3\sqrt n\log n}{\delta_p}\big) = O(\sqrt n)=O(\|E\|).
\]
This proves Theorem \ref{toy1}. The routine verification that  Assumption~{\bf C0} holds in the setting of Theorem \ref{toy1} is deferred to the end of this section.

Before presenting more results, let us pause and make few important remarks. 

\begin{remark} \label{sharpness} (Sharpness) 
As the examples above demonstrate, we 
can often reduce the bound $O(\sigma_{p+1} + \|E \| )$ to $O( \| E \| )$. A natural question arises: {\it can we go even further, reducing $O(\| E \| )$ to something smaller? } 

The answer is no, by results from random matrix theory. Assume that $\frac{m}{n}=\Theta(1)$ and $r=O(1)$. Assume  that $E$ has iid Gaussian $\mathcal{N}(0,1)$ entries and 
$\sigma_1 = \|A\| = c \sqrt n $, for some 
constant $c >2$.  By the reverse triangle inequality, 
\begin{equation} \label{Sharpness0}
    \|\tilde{A}_p -A_p\| \geq | \|\tilde{A}_p \| -\|A_p \| | = |\tilde{\sigma}_1 - \sigma_1|.
\end{equation}

On the other hand, a well-known theorem in random matrix theory, \cite[Theorem 2.8]{BenaychGeorgesNadakuditi2012},  shows that with high probability, 
\begin{equation} \label{sharpness1} \textstyle
    |\tilde{\sigma}_1 - \sigma_1| = \Theta\left( \frac{\|E\|^2}{\sigma_1} \right) = \Theta(\|E\|).
\end{equation}
Therefore, $\| \tilde A_p-A_p\| =\Omega (\| E \|)$. 
\end{remark}

\begin{remark}[A bound via the Davis-Kahan theorem] \label{remark: daviskahan}
Another approach to bounding $\|\tilde A_p-A_p\|$ is to use the Davis--Kahan theorem, which yields a bound of the form
\(
O\left(\frac{\|E\|\sigma_1}{\delta_p}\right)
\); see Appendix \ref{subsec: dktolowrank} for a formal proof. If we want the bound to become $O(\|E \|)$, we need to provide a huge gap condition $\delta_p = \Omega (\sigma_1)$, namely, the gap between the $p$-th and $(p+1)$-th singular value should be of the same order as the largest singular value $\sigma_1$. This is rarely guaranteed, unless the singular values decay exponentially fast.

On the other hand, all of our results, including those in Section~\ref{sec: intro}, allow $\delta_p$ to be significantly smaller than $\sigma_p$. In all of our settings, if we use the Davis-Kahan theorem, the bound obtained would be significantly weaker. For instance, under the assumption of Theorem \ref{toy1} that $\delta_p \ge Cr \log n \cdot \max \big\{ r^2 , \frac{  \sigma_p }{\sqrt{n}}, \frac{n}{\sigma_p} \big\}$ and $\sigma_1=O(\sigma_p)$, Davis-Kahan bound would give 
\[ \textstyle
\|\tilde{A}_p -A_p\|=O\big( \frac{\|E\|\sigma_1}{\delta_p} \big) = O\big( \frac{\|E\|\sigma_p}{\delta_p} \big) = O \big(\min \big\{\frac{\sqrt{n} \sigma_p}{r^3\log n}, \frac{n}{r \log n}, \frac{\sigma_p^2}{r \sqrt{n} \log n} \big\} \big).
\]
In the regime $r=\tilde O(1)$ and $\sigma_p=\Theta(n)$, this gives a bound of order $\tilde O(n)$, whereas our result is of order $\|E\|=\Theta(\sqrt{n})$.

Recent refinements of the Davis--Kahan theorem, such as \cite{OVK13, OVK22, DKTranVu1, JW1}, can further improve the classical bound
\(
O\left(\frac{\|E\|\sigma_1}{\delta_p}\right).
\)
Roughly speaking, under comparable gap assumptions, the resulting bounds remain weaker than Theorem~\ref{theo: mainTh-1} by a factor of $\sigma_1/\sigma_p$. Since a precise comparison requires additional technical assumptions and notation, we omit the details.
\end{remark}

\subsection{Removing \texorpdfstring{$K$}{K}-boundedness} 
The assumption that the entries of $E$ are $K$-bounded is convenient but somewhat restrictive. For entries with unbounded support, this boundedness condition can be enforced with high probability via a simple truncation argument. For example, suppose that the entries $\xi_{ij}$ of $E$ are mean-zero sub-Gaussian random variables with variance at most $\sigma^2$. Then, for any $t>0$,
\( \textstyle
\P(|\xi_{ij}|>t)
\le
2\exp\left(-c\frac{t^2}{\sigma^2}\right)
\)
for some universal constant $c>0$. Hence, by a union bound,
\[ \textstyle
\P\left(
|\xi_{ij}|\le K
\ \text{for all }1\le i\le m, 1 \leq j \leq n
\right)
\ge
1-2n^2\exp\left(-c\frac{K^2}{\sigma^2}\right).
\]
Thus, incurring this additional failure probability, we may apply Theorem~\ref{theo: mainTh-1} with $K$-bounded entries, yielding the following corollary.

\begin{corollary}
Let $A$ be an $m \times n$ matrix of rank $r$  ($m \leq n$) and $E$ be a random matrix, 
whose entries are iid  sub-Gaussian random variables  $\xi_{ij}$ with zero mean, satisfying 
$$\textstyle  \P \left(|\xi_{ij}| > t \right) \leq 2 \exp \left(-c\frac{t^2}{\sigma^2} \right), $$ for some constant $c >0$.  Assume that  \textbf{C0} holds for some parameter $t_1$.  Then for any $t_2, K >0$, we have 
 $$\Norm{\tilde{A}_p - A_p} \leq 32  \sigma_p \left(  \epsilon_1 + rt_1 \epsilon_2 + \frac{r^2 t_2 }{\|E\|} \epsilon_1 \epsilon_2 \right)= 32 \left(\|E\|+ r t_1 \frac{\sigma_p}{\delta_p} + \frac{r^2 t_2 }{\delta_p} \right),$$
\noindent with probability at least $1 - 2 n^2 \exp \left(-c\frac{K^2}{\sigma^2} \right) - r^2 \exp\left( \frac{-t_1^2/2}{\sigma^2 + Kt_1} \right) - \frac{5r(r-1)(m+n) m_4}{t_2^2}. $ 

\end{corollary}

In applications, we can apply this corollary in a manner similar to Theorem \ref{theo: mainTh-1}, by first choosing $K$, and then $t_1$ and $ t_2$ to guarantee a desired probability bound. 

\subsection{Random noise with non-iid entries}
 Now we turn to the case where the entries of $E$ are independent random variables, but not iid.
As the reader will see in the proof, we need to control two types of quantities, $ Ev_k \cdot Ev_l$ and $E^\top u_i \cdot E^\top u_j$, which measure the interaction of the noise $E$ with the singular vectors of $A$.
(Here and later $u \cdot v$ denotes the inner product of the vectors $u$ and $v$.)  When $E$ has iid entries, 
it is easy to show that these quantities have mean zero and are, with high probability, negligible. However, in the general case, their means are not zero. We  introduce  a new parameter $\mu$ to control these means: 
$$\mu:= \max_{1 \leq k<l \leq r} \norm{ \E \left( Ev_k \cdot Ev_l \right)} + \norm{\E \left( E^\top u_k \cdot E^\top u_l \right)} .$$
\noindent A direct computation gives the equivalent representation of $\mu$ as 
$$ \max_{1 \leq k < l \leq r } \bigg\{ \bigg|\sum_{i'=1}^{m} \sum_{i=1}^{n} \sigma_{i' i}^2 v_{li} v_{ki} \bigg|+ \bigg|\sum_{j'=1}^{n} \sum_{j=1}^{m} \sigma_{j j'}^2 u_{lj} u_{kj} \bigg| \bigg\}. $$

 We obtain the following extension of Theorem \ref{theo: mainTh-1}, with an adjustment 
 in the last error term.

\begin{theorem} \label{mainTh01}
Let $A$ be an $m \times n$ matrix of rank $r$  and $E$ be a random matrix, whose entries 
are independent $(K,\sigma)$-bounded  variables with zero mean. 
Assume that \textbf{C0} holds for some parameter $t_1$.  Then for any $t_2 >0$, we have 
$$\Norm{\tilde{A}_p - A_p} \leq 32 \sigma_p \left(  \epsilon_1 + rt_1 \epsilon_2 + \frac{r^2 (t_2+\mu)}{\|E\|} \epsilon_1 \epsilon_2 \right),$$
with probability at least $1 - r^2 \exp\left( \frac{-t_1^2/2}{2\sigma^2 + Kt_1} \right) - \frac{5r(r-1)n m_4}{2t_2^2}.$
\end{theorem}



If we wish to avoid involving the parameter $\mu$, we can use the trivial bound that 
$$\max\{  |Ev_k \cdot Ev_l|,  |E^\top u_i \cdot E^\top u_j| \} \le \| E\|^2. $$ This way, we obtain the following weaker, but more user-friendly, version of Theorem \ref{mainTh01}.
\begin{theorem} \label{cor: AppTrivial}
    Let $A$ be an $m \times n$ matrix of rank $r$ ($m \leq n$) and $E$ be a $(K,\sigma)$-bounded random matrix.
Under Assumption \textbf{C0}, for any $t_1 >0$, we have
$$\Norm{\tilde{A}_p - A_p} \leq 32 \sigma_p \left(\epsilon_1 + r t_1 \epsilon_2 + \frac{r^2 \|E\|^2}{\sigma_p \delta_p} \right) = 32 \left[ \|E\| + \frac{r t_1 \sigma_p }{\delta_p}+ \frac{r^2\|E\|^2}{\delta_p}   \right],$$
 with probability at least $1 - r^2 \exp\left( \frac{-t_1^2/2}{2\sigma^2 + Kt_1} \right).$ 
\end{theorem}

\subsection{Proofs of Theorem \ref{toy1} and Theorem \ref{toy3}} In this subsection, we derive Theorem \ref{toy1} and Theorem \ref{toy3} from Theorem \ref{theo: mainTh-1}. 



 As discussed after Theorem~\ref{theo: mainTh-1}, the error bound in Theorem~\ref{toy1} follows from Theorem~\ref{theo: mainTh-1} by setting
\[
t_1=10K\log n,
\qquad
t_2=r\sqrt n\log n.
\]
It remains to verify Assumption~{\bf C0} under the condition
\begin{equation}\label{toy1 condition} \textstyle
\delta_p
\ge
Cr\log n\cdot
\max \big\{
r^2,\frac{\sigma_p}{\sqrt n},\frac{n}{\sigma_p}
\big\}.
\end{equation}

For the first condition, using $\|E\|=O(\sqrt n)$ and \eqref{toy1 condition},
\[ \textstyle
\frac{\|E\|}{\sigma_p} = O\left(\frac{\sqrt n}{\sigma_p}\right)
=
O\left(\frac{\delta_p}{r\sqrt n\log n}\right)
=o(1).
\]
Similarly, also by $\|E\|=O(\sqrt{n})$ and \eqref{toy1 condition},
\[ \textstyle
r^2t_1\epsilon_2 = O\left(\frac{r^2\log n}{\delta_p}\right)
=o(1),\,\text{and}\,\,\sqrt{r\eta} =\sqrt{\frac{r\|E\|^2}{\sigma_p\delta_p}}
=
O\left(
\sqrt{\frac{rn}{\sigma_p\delta_p}}
\right)
=o(1).
\]
 Hence, Assumption~{\bf C0} holds, completing the proof of Theorem~\ref{toy1}.

Next, we prove Theorem \ref{toy3} -  another analogue of Theorem \ref{theo: mainTh-1} with $t_1=10K\log n$ and $
t_2=r\sqrt n\log n.$ First, the assumption of Theorem \ref{toy3} that $r \leq \log^{-1} n \cdot \min\{\sqrt{\delta_p}, \frac{\delta_p \sigma_p}{n}  \} $ implies 
\[\textstyle
r^2 t_1 \epsilon_2 = O \big( \frac{r^2 \log n}{\delta_p} \big) = O (\log^{-1} n) =o(1),\,\text{and}\,\,\sqrt{r\eta}=O\big(\sqrt{\frac{r n}{\sigma_p \delta_p}} \big)=O(\sqrt{\log^{-1}n})=o(1).
\]
Thus Assumption \textbf{C0} is satisfied. 

For the error bound of Theorem \ref{toy3}, using $\sigma_p=\sigma_{p+1}+\delta_p$ and $\sigma_p \geq \sqrt{r^{-1} mn}$, the RHS of Theorem \ref{theo: mainTh-1} becomes
\begin{equation*}
    \begin{split}
        32 \bigg(\|E\|+ rt_1 \frac{\sigma_p}{\delta_p} + \frac{r^2 t_2}{\delta_p} \bigg) & = O \bigg(\|E\| + rt_1 \frac{\sigma_{p+1}+\delta_p}{\delta_p} + r^2 \frac{t_2}{\delta_p} \bigg) \\
        & = O \bigg[\|E\| + \sigma_{p+1} \cdot \big( \frac{r \log n}{\delta_p} + \frac{r^3 \sqrt{n} \log n}{\sigma_{p} \delta_p} \big)  \bigg] \\
        & = O \bigg[\|E\| + \sigma_{p+1} \cdot \big( \frac{r \log n}{\delta_p} + \frac{r^{7/2} \log n}{\sqrt{m} \delta_p} \big)  \bigg] = O \big(\|E\|+\epsilon_0 \sigma_{p+1} \big).
    \end{split}
\end{equation*}
This proves Theorem \ref{toy3}.


\section{Deterministic results } \label{subsec: mainDe}

\subsection {A deterministic theorem}

In this section, we present a result which applies to a deterministic $E$. The random-noise results are obtained from this deterministic theorem using estimates from random matrix theory; see Appendix \ref{xy random}.

We first introduce two key parameters,  $x$ and $y$, which account for the interaction between the singular vectors of $A$ and $E$. Here and later, $u \cdot v$ is the inner product of two vectors $u$ and $v$. 
\begin{definition} \label{def: xyz}
\end{definition}
\begin{itemize}
\item $x:= \max_{1 \leq i,j \leq r} |u_i^\top E v_j|$.
\item $y:= \max_{1\leq i < j\leq r} \left\lbrace \norm{Ev_i \cdot Ev_j},\norm{ E^\top u_i \cdot E^\top u_j} \right\rbrace $.
\end{itemize}

The parameters measure the intensity of $E$ with respect to the singular vectors of $A$.

\begin{theorem} \label{cor: rec} Let $A$ be an $m \times n$ matrix of rank $r$ and let $p \leq r$ be a natural number. If 
\begin{equation} \label{ass: cormain}
   \max \left\lbrace  \frac{p\|E\|}{\sigma_p}, \frac{r^2 x}{\delta_p} , \frac{ \sqrt{r} \|E\|}{\sqrt { \sigma_p \delta_p }} \right\rbrace   \le \frac{1}{96 } , 
\end{equation}
  then  
  \begin{equation} \label{finalboundHRApp-friendly2}
      \big\|\tilde{A}_p -A_p\big\| \leq 32 \sigma_p \left(\frac{p\|E\|}{\sigma_p} + \frac{pr x}{\delta_p} + \frac{r^2 y}{\sigma_p \delta_p} \right).
  \end{equation}
\end{theorem}

The constants $32$ and $96$ are ad hoc, and both can be reduced by optimizing the chosen constants in Section~\ref{sec: proof}. If we replace \eqref{ass: cormain} by  a stronger assumption that 
\begin{equation} 
   \max \left\lbrace  \frac{p\|E\|}{\sigma_p}, \frac{r^2x}{\delta_p} , \frac{ \sqrt{r} \|E\|}{\sqrt { \sigma_p \delta_p }} \right\rbrace \le \epsilon, 
\end{equation} where $\epsilon >0$ is a small constant, then the bound becomes 
\begin{equation} \label{finalboundHRApp-friendly2-new}
      \big\|\tilde{A}_p -A_p\big\| \leq (6+O(\epsilon))\sigma_p \left(\frac{p\|E\|}{\sigma_p} + \frac{pr x}{\delta_p} + \frac{r^2 y}{\sigma_p \delta_p} \right).
  \end{equation}



\subsection{   Weakening the low-rank assumption} \label{subsec: notlowrank}
We can extend Theorem \ref{cor: rec} to th full rank setting  by redefining the parameter $r$. 
In \cite{DKTranVu2}, we prove a general theorem which 
allows the following choice for $r$ (instead of being the rank). 
 Given $p$,  we define $r(p)$ as the largest integer such that 
$$ \sigma_{r(p)} \ge \frac{1}{2} \sigma_p. $$

 Theorem \ref{cor: rec} (and all other new results in this paper) holds when we let $r$ be this ``halving" parameter $r(p)$ instead of the rank of $A$. The constant $1/2$ can be replaced by any positive constant less than 1. 
 In order to make $r(p)$ relatively small, we only need to assume that the singular values after $\sigma_p$ decay reasonably fast.   In the opposite case, if $r(p)$ is large, then there are too many singular values between $\sigma_p$ and $\sigma_p/2$. This implies that $A_p$ may not be a good approximation of $A$, as there is too much energy 
 of the matrix left beyond the first $p$ dimensions. 

 The proof of this stronger result, however,  is substantially more technical than that of Theorem \ref{cor: rec}. We will present the detailed result and proof in a future paper \cite{DKTranVu2}. In that paper, we will analyze the perturbation of a general function $f(A)$. The case $f(z)=z$ corresponds to the low-rank perturbation. 
 Readers interested in \cite{DKTranVu2}  may find the relatively shorter proof presented here useful as a warm-up. Related arguments for the special case $f \equiv 1$ appeared in \cite{DKTranVu1}.

\subsection {The  parameters \texorpdfstring{$x$}{x} and \texorpdfstring{$y$}{y}: the blessing of high dimensions } \label{section:skewness}
A central idea in this paper and our related work  \cite{DKTranVu1, DKTranVu2} is to incorporate the (skewness) parameters $x,y$, which measure the relative strength of the noise matrix $E$ with respect to the singular vectors of $A$, as part of perturbation bounds instead of using the absolute intensity $\|E\|$. Notice that one can trivially bound $x$ and $y$ by $\|E\|$ and $\|E\|^2$,  respectively. However, in the random setting, it is very easy to significantly improve these bounds. 

When $E$ is random and $u$ is a fixed vector, the vector $Eu$ behaves like a random vector. The point is that in high dimensions, a random vector and a fixed vector, or two random vectors, are very close to being orthogonal. We can exploit this ``skewness" to obtain a good bound on the inner product of these vectors,  often gaining a factor $n^{1/2-o(1)}$, compared to the trivial bound, where $n$ denotes the ambient dimension.
Thus, the higher the dimension, the more we gain. In other words, in this situation, the high dimension is not a curse, but a blessing.

 {\it \noindent Example.} The reader can prove the following as an exercise, or see \cite{Ver1book}. Assume that 
 $E$ is an $n \times n$ random matrix with iid sub-Gaussian entries with mean zero and variance one, and $u$ and $v$ are two orthogonal unit vectors in $\mathbb{R}^n$, and $m = \Theta(n)$. Then 
 with probability $1-o(1)$

 \begin{itemize} 

 \item $\| E \| = \Theta (\sqrt n) $.

 \item  $ \| E u \|, \|Ev\|  =\Theta (\sqrt n)$. 
 
 \item  $| u^\top E u |$, $| u^\top E v | = O(\log n)$. 
 
 \item  $ |Eu \cdot Ev | = O(\sqrt n \log n )$.

\item $\log n $ can be 
reduced to $\omega(1 )$ for a fixed pair $u$ and $v$. We use
$\log n$ for convenience, taking into account that the bound still holds when taking maximum over polynomially many choices of $u$ and $v$.

 \end{itemize}

Finally, let us point out that even with the 
trivial bounds ($\| E \|$ for $x$ and $\|E\|^2$ for $y$),  
Theorem \ref{cor: rec} still provides useful information. These trivial estimates imply 
$$\textstyle \big\|\tilde{A}_p -A_p\big\| \leq 32 \sigma_p \left(\frac{p\|E\|}{\sigma_p} + \frac{pr \|E\|}{\delta_p} + \frac{r^2 \|E\|^2}{\sigma_p \delta_p} \right) = O \left(\frac{r^2 \sigma_p \|E\|}{\delta_p} \right). $$
Since $\sigma_p =\sigma_{p+1}+\delta_p$, we can rewrite this inequality as 
$$\textstyle\big\|\tilde{A}_p -A_p\big\| = O \left(r^2 \|E\| + \sigma_{p+1} \frac{r^2 \|E\|}{\delta_p}\right). $$
 This simple corollary   still improves upon  Theorem \ref{EY0}  in certain regimes.

\subsection{  Small gap assumptions} \label{gap-nec}
We next discuss the gap assumption on $\delta_p$. 
First, we would like to point out that {\it some} gap assumption is necessary.  If $\sigma_p $ and $\sigma_{p+1}$ are too close, then there is a non-negligible probability that their order is reserved after perturbation. This would cause  $\| \tilde A_p - A_p \| $ to be large.  
 
 For instance, suppose that these two singular values are swapped, while the remaining and singular vectors change only negligibly, then   $ \tilde A_p - A_p $ is essentially 
$\sigma_p u_p v_p^\top - \tilde \sigma_p u_{p+1} v_{p+1}^\top $ (the crucial point is  the appearance of $u_{p+1} $ and $v_{p+1} $ after the swap). On the other hand, 
\begin{equation} \label{EY2}  \| \sigma_p u_p v_p^\top - \tilde \sigma_p u_{p+1} v_{p+1}^\top  \| = \max \{\sigma_p, \tilde\sigma_{p}  \} \approx \max \{\sigma_p, \sigma_{p+1} +\|E\|\}. \end{equation}

The question 
here is how {\it strong} the gap assumption must be. 
In existing analyses concerning a random $E$, the common assumption (after some normalization if necessary) is that 
 \begin{equation} \label{strongassumption} \delta_p \gg \|E\| = \Theta (\sqrt n); \end{equation} e.g., one needs this to make a non-trivial use of the Davis-Kahan bound; see Theorem \ref{DKtheorem}. On the other hand, in many typical regimes (e.g., the setting in Subsection \ref{newsimplebound}), 
 our assumption \eqref{ass: cormain} only requires 
 \begin{equation} \label{weakassumption} \delta_p \ge C r^2\log n, \end{equation}  which is superior for $r =o( n^{1/4})$.
 In particular, if $r $ is polylogarithmic (in $n$), then we only need to assume $\delta_p$ to be polylogartihmic as well.

\subsection{A geometric viewpoint  } \label{subsec: geoview}
 It is interesting to view the new result geometrically, with respect to the triangle inequality used in the proof of Theorem \ref{EY0}.  We view a matrix as a point in 
 $R^{m \times n}$. 

\usetikzlibrary{decorations.pathreplacing}
$$\begin{tikzpicture}
\coordinate (A0) at (1,0);
\node[below] at (A0){$\tilde{A}_p$};
\coordinate (A) at (6.75,2);
\node[below] at (A){$A_p$};
\coordinate (A1) at (2,3);
\node[above] at (A1){$\tilde{A}$};
\coordinate (A2) at (3.25,3);
\node[above] at (A2){$A$};
\coordinate (A3) at (9.5,0);
\node[below] at (A3){$\tilde{A}_p$};
\coordinate (A4) at (11.25,0.25);
\node[below] at (A4){$A_p$};
\coordinate (A5) at (8,3);
\node[above] at (A5){$\tilde{A}$};
\coordinate (A6) at (9.25,3);
\node[above] at (A6){$A$};
\coordinate (B) at (10.15,0.25);
\node[above] at (B){$O(\|E\|)$};
\coordinate (B1) at (8.75, 3);
\node[below] at (B1){$\|E\|$};
\coordinate (B2) at (2.65, 3);
\node[below] at (B2){$\|E\|$};
\coordinate (C) at (10.5, 1.75);
\node[right] at (C){$\sigma_{p+1}$};
\coordinate (C1) at (5.15, 2.65);
\node[right] at (C1){$\sigma_{p+1}$};
\coordinate (D) at (8.75, 1.5);
\node[left] at (D){$\tilde{\sigma}_{p+1}$};
\coordinate (D1) at (1.5, 1.75);
\node[left] at (D1){$\tilde{\sigma}_{p+1}$};

\draw[very thick, red] (A0) -- (A1);
\draw[very thick, red] (A) -- (A2);
\draw[very thick, red] (A3) -- (A5);
\draw[very thick, red] (A4) -- (A6);
\draw[very thick, blue, dashed] (A0) -- (A);
\draw[very thick, blue, dashed] (A1) -- (A2);
\draw[very thick, blue,dashed] (A3) -- (A4);
\draw[very thick, blue, dashed] (A5) -- (A6);

\filldraw (A1) circle (2pt);
\filldraw (A0) circle (2pt);
\filldraw (A) circle (2pt);
\filldraw (A3) circle (2pt);
\filldraw (A2) circle (2pt);
\filldraw (A4) circle (2pt);
\filldraw (A5) circle (2pt);
\filldraw (A6) circle (2pt);


\end{tikzpicture}.$$

We assume that the distance $\| E \|$ between $A$ and $\tilde A$ is relatively small compared to the distance $\sigma_{p+1}$ between $A$ and $A_p$. Since
$$| \tilde \sigma_{p+1} - \sigma_{p+1} | \le \| E\| \ll \sigma_{p+1},$$ 
we conclude that the lengths of the two segments $AA_p$ and $\tilde A \tilde A_p$ are roughly the same, and are much longer than the length of the segment $A \tilde A$.

It is natural to think that the segments $AA_p$ and 
$\tilde A  \tilde A_p $ are in general position with respect to each other. In this case, it is easy to see that the segment $A_p \tilde A_p$ is also of order $\Theta (\sigma_{p+1})$, which is much larger than $A\tilde A$. This corresponds to the configuration shown in the left panel. 

\vskip2mm 

\noindent{\it Exercise.} Assume that $ \| A\tilde A \|  =\epsilon =o(1)$, $ \| AA_p\| = \| \tilde A \tilde A_p \| =1$. Fix the points $A, A_p, \tilde A$ and let $\tilde A_p$ be a random point on the unit sphere centered at $\tilde A$. Then with high probability, $\| A_p \tilde A_p \| = (1+o(1)) \sqrt 2$.

\vskip2mm

The right panel illustrates the behavior predicted by our findings. 
In the cases when the term $\sigma_px/\delta_p$  is small compared to $\| E\|$, our bound is 
$\|  \tilde A_p - A_p \| = O(\| E \| )$. This means that the 
length of the segment $A_p \tilde A_p$ is of order $\| E \|$. This implies that the two segments  $AA_p$ and $\tilde A \tilde A_p$  {\it are not} in general position with respect to each other, as one might initially believe. Instead, \textit{they must be nearly parallel}.

\section {The symmetric version } \label{section: symmetric}
In this section, we discuss the case when $A$ is a symmetric matrix of rank $r$. We can reduce the general case when $A$ is an $m \times n$ matrix to this case by a standard symmetrization argument. Furthermore, our proof technique 
makes essential use of symmetry. 

When $A$ is symmetric,  it admits the spectral decomposition $A = \sum_{i=1}^r \lambda_i u_i u_i^\top$, where $\lambda_i$ are the eigenvalues and $u_i$ are the (corresponding) eigenvectors.  The (nontrivial) singular values are $|\lambda_i |$ with
$i=1, 2, \dots, r$.

For simplicity, we first consider the case when $A$ is positive semi-definite. In this case, there  is no difference between singular values and eigenvalues, as  $\sigma_i = \lambda_i$ for all $1 \leq i \leq r$,  and 
$$A_p = \sum_{i=1}^p \sigma_i u_i u_i^\top = \sum_{i=1}^p \lambda_i u_i u_i^T.$$ 

\begin{theorem} \label{theo: mainPositive} Let $A$ be a positive semi-definite  matrix of rank $r$ and let $p \leq r$ be a natural number. If 
\begin{equation} \label{ass: sym0}
 \max \left\lbrace  \frac{p\|E\|}{\sigma_p}, \frac{r^2x}{\delta_p} , \frac{ \sqrt r \|E\|}{\sqrt { \sigma_p \delta_p }} \right\rbrace   \le \frac{1}{24} ,   
\end{equation}
  then  
$$ \big\|\tilde{A}_p -A_p\big\| \leq 8 \sigma_p \left(\frac{p\|E\|}{\sigma_p} + \frac{pr x}{\delta_p} + \frac{r^2 y}{\sigma_p \delta_p} \right). $$
\end{theorem}

\begin{remark} Theorem \ref{cor: rec} involves  larger constants ($96,32$) than Theorem \ref{theo: mainPositive} ($24,8$), 
due to the symmetrization step; see Subsection \ref{subsec: proofrec} for details.
\end{remark}
If we replace \eqref{ass: sym0} by  a stronger assumption that 
\begin{equation} \textstyle
   \max \left\lbrace  \frac{p\|E\|}{\sigma_p}, \frac{r^2x}{\delta_p} , \frac{ \sqrt{r} \|E\|}{\sqrt { \sigma_p \delta_p }} \right\rbrace \le \epsilon, 
\end{equation} where $\epsilon >0$ is a small constant, then the bound becomes 
\begin{equation} \label{finalboundHRApp-friendly2-new-sym}
      \big\|\tilde{A}_p -A_p\big\| \leq (3+O(\epsilon))\sigma_p \left(\frac{p\|E\|}{\sigma_p} + \frac{pr x}{\delta_p} + \frac{r^2 y}{\sigma_p \delta_p} \right).
  \end{equation}

Now we turn to the general case in which $A$ has both positive and negative eigenvalues. 
We identify the leading $p$ singular values as follows. 
There is an integer  $0 \le k \le p$ such that
$$ \lambda_1 \geq \lambda_2 \geq \cdots \geq \lambda_k > 0 \geq \lambda_{n-(p-k)+1} \geq \cdots \geq \lambda_n ,$$
and for $S:= \{1, \dots, k, n- (p-k) +1, \dots, n \} $, we have
$$\{ \sigma_1 , \dots, \sigma_p \}  = \{ | \lambda_i |, i \in S \} . $$ 
Thus, instead of using the gap  $\delta_p:= \sigma_p- \sigma_{p+1}$, now we need to define a new eigenvalue gap with respect to $S$:
\begin{equation} \label{gapdef}
      \delta_S :=\min\{\lambda_k -\lambda_{k+1}, \lambda_{n-(p-k)} - \lambda_{n-(p-k)+1}\}.
 \end{equation}

\begin{theorem} \label{theo: main} Let $A$ be a symmetric matrix of rank $r$ and let $p \leq r$ be a natural number. If 
\begin{equation*} 
 \max \left\lbrace  \frac{p\|E\|}{\sigma_p}, \frac{r^2 x}{\delta_p} , \frac{ \sqrt r \|E\|}{\sqrt { \sigma_p \delta_p }} \right\rbrace   \le \frac{1}{24} ,   
\end{equation*}
then
  \begin{equation} \label{finalboundHRApp-friendly2-theo}
      \big\|\tilde{A}_p -A_p\big\| \leq 8 \sigma_p \left(\frac{p\|E\|}{\sigma_p} + \frac{pr x}{\delta_S} + \frac{r^2 y}{\sigma_p \delta_S} \right).
  \end{equation}
\end{theorem}
Note that the gap $\delta_S$ is at least $\delta_p$, but it can be significantly bigger because $A$ might have a small singular value gap even when all relevant eigenvalue gaps are large. 


\section{Proofs of the deterministic theorems: General strategy} \label{sec: proofidea}
In this section, we first outline the key ideas behind the proof of Theorem~\ref{theo: mainPositive}. 
The proof of Theorem~\ref{theo: main} then follows with a straightforward modification; see Remark~\ref{modification}. 
Finally, Theorem~\ref{cor: rec} is derived from Theorem~\ref{theo: main} using a standard symmetrization argument. 

Before turning to the main proofs, we recall Cauchy's integral theorem~\cite{CAbook}, which will be used in several of our estimates.  
\begin{theorem}[Cauchy's integral theorem] \label{theo: Cauchy}
Let $\Gamma$ be a simple closed contour, and let $f$ be an analytic on a simply connected domain $S$ containing $\Gamma$ and its interior. Then,
$$
   \frac{1}{2 \pi {\bf i}} \int_{\Gamma} \frac{f(z)}{zI-A}\,dz 
   = 
   \begin{cases}
      f(a), & a \text{ inside } \Gamma, \\[4pt]
      0, & a \text{ outside } \Gamma.
   \end{cases}
$$
\end{theorem}

\subsection{Proof of Theorem \ref{theo: mainPositive}: A sketch}  


In the setting of this theorem,  $A$ is a positive semi-definite symmetric matrix of rank $r$, which admits the spectral decomposition 
$$A= \sum_{i=1}^n \lambda_i u_i u_i^\top.$$
Here, 
$\lambda_1 \geq \lambda_2 \geq \cdots \geq \lambda_p > \lambda_{p+1} \geq \cdots \geq \lambda_n \geq 0$ are the eigenvalues and $u_1, u_2, \dots, u_n$ are the corresponding eigenvectors. The best rank-$p$ approximation of $A$ is 
$$A_p:= \sum_{i=1}^p \lambda_i u_i u_i^\top.$$
Similarly, we have 
$\tilde{A}=\sum_{i=1}^n \tilde{\lambda}_i \tilde{u}_i \tilde{u}_i^\top,$ with $\tilde{\lambda}_1 \geq \tilde{\lambda}_2 \geq \cdots \geq \tilde{\lambda}_n$. Let $\tilde{A}_p$ be the best rank-$p$ approximation of $\tilde{A}$.  Note, however, that 
$\tilde{A}$ is not necessarily positive semi-definite.

We are going to adapt the contour approach from \cite{DKTranVu1}, with a significant modification that exploits the low rank assumption. Applying Cauchy's integral theorem (Theorem \ref{theo: Cauchy}) for a simple closed contour $\Gamma$ and for $f(z)=z$, we have  
$$\frac{1}{2 \pi \mathrm{i}} \int_{\Gamma} \frac{z}{zI-A}\,dz 
   = \begin{cases}
      a, & a \text{ inside } \Gamma, \\[4pt]
      0, & a \text{ outside } \Gamma.
   \end{cases}$$
 
 Let $\Gamma$ be a contour enclosing eigenvalues $\{\lambda_i : 1 \le i \le p \}$, with all $\lambda_j$ for $j > p$ outside $\Gamma$. Then by the identity above (see also \cite{Book1,Kato1}), we have 
\begin{equation} \label{contour-formula} 
  \frac{1}{2 \pi {\bf i}} \int_{\Gamma} z (z I-A)^{-1}\,dz 
  = \sum_{i=1}^p  \lambda_i u_i u_i^\top.
\end{equation}

To be specific, we can  construct $\Gamma$ as a rectangle with  
\begin{itemize}
    \item horizontal edges parallel to the real axis,
    \item the left vertical edge bisecting the interval $(\lambda_{p+1}, \lambda_p)$,
    \item the right vertical edge passing through $\|A\| + 1.1\,\|E\|$.
\end{itemize}

Now, let $\tilde{S} \subset [n]$ denotes the set of indices $i$ for which $\tilde \lambda_i$ lies inside $\Gamma$, and all $\tilde \lambda_j, j \notin \tilde{S}$ are outside.\footnote{If some eigenvalues of $\tilde A$ lie on the left vertical edge of $\Gamma$, we may shift this edge slightly so that it passes through
\(
(1+10^{-c})\frac{\lambda_p+\lambda_{p+1}}{2}
\)
instead of exactly bisecting the interval $(\lambda_p,\lambda_{p+1})$.} Define $\tilde{\Pi}_{\tilde{S}}$ as the orthogonal projection onto the space spanned by eigenvectors $\{ \tilde{u}_i, i \in \tilde{S}\}$. Under Assumption \eqref{ass: sym0} of Theorem \ref{theo: mainPositive}, \cite[Theorem 2.4]{DKTranVu1} applies, yielding
$$\| \tilde{\Pi}_{\tilde{S}} - \Pi_p\| \leq  12 \sqrt{p} \left( \frac{ \|E\|}{\lambda_p} + \frac{ \sqrt{r} x}{\delta_p} + \frac{\sqrt{r} \|E\|^2}{\delta_p \lambda_p} \right) < \frac{1}{3 \sqrt{p}} + \frac{1}{3 r} + \frac{1}{48} < 1.$$
This enforces $\tilde{\Pi}_{\tilde{S}} $ and $ \Pi_p$ to have the same rank, i.e., $|\tilde{S}| = p$. Moreover, since $\Gamma$ extends beyond $\|A\| +1.1  \|E\|$, the Weyl's inequality (Theorem \ref{stanTheo: Weyl}) implies that $\tilde{\lambda}_1$. Therefore, $\tilde{S}$ must be $ \{1,2, \dots, p\}$. Equivalently, $\Gamma$ encloses exactly the  $p$ leading eigenvalues of $\tilde{A}$. 

Furthermore, by the construction of $\Gamma$ and the
assumption \eqref{ass: sym0}, we have 
$$\tilde{\lambda}_{p} \ge  \frac{\lambda_{p+1}+\lambda_p}{2}\ge \lambda_p/2 \ge 2 \| E \|. $$ 
On the other hand, by Weyl's inequality and the fact that $A$ is positive semi-definite, 
$$\tilde \lambda_n \ge \lambda_n - \| E\| \ge - \| E \|.  $$
Since $\tilde{\lambda}_p \geq 2\|E\|$ while every negative eigenvalue has magnitude at most $\|E\|$, the first $p$ eigenvalues of $\tilde A$ are precisely the $p$ leading singular values. Thus, 
$\Gamma$ exactly encloses the $p$ leading singular values of $\tilde{A}$. Therefore, 
 \begin{equation} \label{contour-formula1} 
  \tilde A_{p} = \sum_{ i \in \tilde{S}} \tilde \lambda_i \tilde u_i  \tilde u_i ^\top =\frac{1}{2 \pi {\bf i}}  \int_{\Gamma}  z (z I-\tilde A)^{-1} dz.  \end{equation}

 \noindent We therefore obtain a key identity 
  \begin{equation} \label{key perturb contour repr}
 \tilde A_{p}- A_p=  \frac{1} {2 \pi {\bf i} } \int_{\Gamma}  z [(z I-\tilde A)^{-1}- (z I- A)^{-1} ]  dz. 
  \end{equation} 

Next, using  the resolvent formula \cite[Chapter 2]{CAbook} repeatedly, we obtain 
\begin{equation} \label{main expansion} (z I-\tilde A)^{-1}-(z I-A)^{-1}   =\sum_{s=1}^{\infty}  (z I-A)^{-1} [ E(z I-A)^{-1} ]^s ,
\end{equation}  where the assumption \eqref{ass: sym0} guarantees that the RHS indeed converges (see \cite[Section 9.5]{DKTranVu1}). Therefore, 
$$ \tilde A_{p}- A_p=  \sum_{s=1}^{\infty} F_s,$$
in which
\begin{equation*}
F_s := \frac{1} {2 \pi {\bf i} }\int_{\Gamma} z (zI-A)^{-1} [ E (zI-A)^{-1} ]^s dz.
\end{equation*}
By the triangle inequality, we further have
\begin{equation} \label{bound00} 
\big\|\tilde A_p- A_p\big\| \le \sum_{s=1}^{\infty} \Norm{F_s}.
\end{equation}

This contour approach has been used by many researchers; see, for instance, Kato's classical book  \cite[Chapter 2]{Kato1}. 
The key task is to sharply bound $\| F_s \|$. In \cite[Chapter 2]{Kato1}, 
one  can use
\begin{equation}\label{trivialF_s}
    \| F_s \|  \le \frac{1} {2 \pi } \|E\|^{s} \max_{z \in \Gamma}|z| \int_{\Gamma} \|(zI-A)^{-1}\|^{s+1} |dz| = O \left[ \max_{z \in \Gamma} |z| \left( \frac{\|E\|}{\delta_p} \right)^s \right]. 
\end{equation}
If $\delta_p > \|E\|$, then the RHS of \eqref{bound00} is summable and we   obtain 
$$ \textstyle \sum_{s=2}^{\infty} \|F_s\| = O \bigg[\max_{z \in \Gamma} |z| \left(\frac{\|E\|}{\delta_p} \right)^2   \bigg] =
O\bigg[ \sigma_1 \left( \frac{\|E\|}{\delta_p} \right)^2 \bigg] .$$
One can show that  $\| F_1 \|$ is negligible compared to $\big[ \sigma_1 \big( \frac{\|E\|}{\delta_p} \big)^2 \big]$ (e.g., see Subsection \ref{subsec: F1}). Thus, we have 
$$\textstyle \big\|\tilde{A}_p -A_p\big\| \leq \|F_1\| + O \left[ \sigma_1 \left( \frac{\|E\|}{\delta_p} \right)^2  \right] = O \left[ \sigma_1 \left( \frac{\|E\|}{\delta_p} \right)^2  \right].$$

This bound is nontrivial and can improve upon Theorem \ref{EY0} in many regimes. 
However, in many settings, it is larger than our new bound (in Theorem \ref{theo: mainTh-1} or in Theorem \ref{theo: mainPositive})  
by a factor $\sqrt n$. Gaining this extra factor of $\sqrt n $ is critical in many applications. 

In \cite{tran2025davis}, assuming that the gap $\delta_p \ge 4 \| E \| $, the current authors developed another approach to bound the RHS of \eqref{key perturb contour repr} using a bootstrapping argument.  To be precise, this contour argument was applied to bound $\| \tilde \Pi_p - \Pi_p \| $, thus there is no factor $z $ on the RHS. In \cite{TranVishnoiVu2025}, 
it is shown that one can use the same trick to handle the general setting where 
$z$ is replaced by any function $f(z)$. With respect to $\| \tilde A_p -A_p \|$, we obtain a bound comparable to the bounds in this paper, up to a logarithmic term. The bootstrapping approach is simple to use, but it depends strongly on the gap assumption $\delta_p \ge 4  \| E\|$. 
The limitation of this bootstrapping analysis lies in bounding 
\begin{equation*}
\| F_s  \| := \big\| \frac{1} {2 \pi {\bf i} }\int_{\Gamma} z (zI-A)^{-1} [ E (zI-A)^{-1} ]^s dz \big \|\,\,\text{by}\,\,\,\frac{1} {2 \pi  }\int_{\Gamma} \big \| z (zI-A)^{-1} [ E (zI-A)^{-1} ]^s \big \|\, |dz|.
\end{equation*} 



We design a new method to estimate  $\| F_s \|$, which involves two extra expansions.
Notice that $ (zI-A)^{-1} =\sum_{i =1}^n  \frac{u_i u_i^\top}{ z-\lambda_i}  $. 
In order to control $F_s$, we use the fact that $A$ is positive semi-definite of rank $r$ and split the RHS  into two parts 
$P= \sum_{1 \le i \le r } \frac{1}{ z -\lambda_i} u_iu_i^\top $ and $Q=  \frac{\sum_{ j > r } u_j u_j^\top} {z}$.

We now expand this product as the sum of $2^{s+1} $ terms, each of which is the form $z R_1 E R_2 \dots, E R_{s+1} $,  where $R_i$ is either $P$ or $Q$. This is the second expansion in our proof, beyond the Taylor expansion used earlier. 

 The only term that does not involve any nontrivial eigenvalue is $z QE \dots EQ $ (which involves only $Q$-operators). Using the formula for $Q$, we further expand this 
product as a sum of $(n- r)^{s+1} $ terms, each  of the form 
$$ \frac{z \times u_{l_1} u_{l_1} ^\top E  u_{l_2} u_{l_2} ^\top E \cdots E u_{l_s} u_{l_s} ^\top}{ z^{s+1} }=\frac{u_{l_1} u_{l_1} ^\top E  u_{l_2} u_{l_2} ^\top E \cdots E u_{l_s} u_{l_s} ^\top}{ z^{s} }   . $$
The key point here is that by the construction of $\Gamma$, the origin is {\it outside } the contour, thus the integral of $    \frac{1}{ z^{s}  }$ is zero for any $s \geq 1$. Therefore, $ \frac{1} {2 \pi {\bf i} }\int_{\Gamma} z QE \dots EQ dz = 0$. 

 Next, we observe that all remaining terms contain many $P$-operators, each of which is adjacent to an $E$. An expansion of the form above 
 will give rise to the terms $x$ and $y$ in Theorem \ref{cor: rec}, which is the main source of our improvement.

For a concrete example, let us consider the term $ z PEP$ ($s=1$). Writing $P= \sum_{1 \leq i \leq r} \frac {1}{ z-\lambda_i } u_i u_i^\top $ and expanding the product (this is the third expansion in our analysis), we obtain 
$$z PEP= \sum_{ 1\le i_1, i_2 \le r  }  \frac{z} { (z-\lambda_{i_1} ) (z -\lambda_{i_2} ) }  u_{i_1} u_{i_1} ^\top  E u_{i_2} u_{i_2} ^ T =  \sum_{1 \le i_1, i_2 \le r }  \frac{z} { (z-\lambda_{i_1} ) (z -\lambda_{i_2} ) } 
 u_{i_1} u_{i_2} ^\top  (u_1^\top  E u_{i_2}) . $$

%
As $u_{i_1}, u_{i_2} $ are unit vectors, a  trivial bound for $|u_{i_1} ^\top E u_{i_2}  |$ is $\| E \| $. However, as mentioned before,  if $E$ is in general position with respect to 
both vectors $u_{i_1} $ and $u_{i_2} $, then we expect a saving.  (If $E$ is random, in particular, then we often achieve 
the $\sqrt n$ saving here.) This naturally leads to the definition of 
the critical ``skewness"  parameter $x$ (see Definition \ref{def: xyz}). 

The main technical part of the proof is to handle $  z (P+Q) [E (P+Q)] ^s $, for a general $s$. After the second and third expansions, it becomes an (enormous) sum of terms of the form 
$$ h(z) u_{l_1} u_{l_1}^\top E u_{l_2} u_{l_2} ^\top E \dots E u_{l_{s+1} } u_{l_{s+1} }^\top,\,\,\text{where}\,\,h(z)= \frac{z}{ \prod_{i_l \in I}  (z- \lambda_{i_l}) \times z^{s+1 - | I| }},$$  
for some sequence $I=[i_1, i_2, \dots]$ of indices in 
$\{1,2, \dots, r \}$ of length at most $s+1$. 


 In order to estimate these terms and sum them, 
 we introduce a key lemma, Lemma \ref{key contour1}, which we prove using induction.
The fact that the ground-truth matrix $A$ has low rank plays a key role in our analysis.

 The contour approach introduced here can be extended to cover the full rank case \cite{DKTranVu1, DKTranVu2}. However, 
 the analysis becomes substantially more complicated and very long. In particular, it requires a delicate combinatorial profiling of the terms by an auxiliary graph. On the other hand, it leads to the extension 
 (with the new definition of $r$) discussed in Subsection \ref{subsec: notlowrank}. The proof presented here could serve as a warm-up to those more technical proofs.

\subsection{Proof of Theorem \ref{theo: main}} \label{modification} 

The proof of Theorem~\ref{theo: main} follows the same overall strategy, except for the construction of the contour $\Gamma$, which must be adapted to guarantee the identities \eqref{contour-formula1} and \eqref{key perturb contour repr}. Recall that there exists an integer $0 \leq k \leq p$ such that
\[
 \lambda_1 \geq \lambda_2 \geq \cdots \geq \lambda_k > 0 \geq 
 \lambda_{n-(p-k)+1} \geq \cdots \geq \lambda_n,
\]
and for 
\[
S := \{1, \dots, k,\, n-(p-k)+1, \dots, n\},
\]
we have
\[
\{\sigma_1, \dots, \sigma_p\} = \{ |\lambda_i| : i \in S \}.
\]
Hence, the key requirement for $\Gamma$ is that it enclose exactly the eigenvalues $\lambda_i, \tilde{\lambda}_i$ with $i \in S$.  

To achieve this, we construct $\Gamma$ as the union of two disjoint rectangles $\Gamma_1 \cup \Gamma_2$, where:  
\begin{itemize}
    \item $\Gamma_1$ isolates $\{\lambda_1, \lambda_2, \dots, \lambda_k\}$, extending to the right beyond $\lambda_1 + \|E\|$, with its left vertical edge intersecting the real axis at $\sigma_p - \delta_p/2$;  
    \item $\Gamma_2$ isolates $\{\lambda_n, \lambda_{n-1}, \dots, \lambda_{n-(p-k)+1}\}$, extending to the left beyond $\lambda_n - \|E\|$, with its right vertical edge intersecting the real axis at $-\sigma_p + \delta_p/2$.  
\end{itemize}
See the figure below for an illustration.

\resizebox{\linewidth}{!}{
\begin{tikzpicture}
\coordinate (A) at (7,0);
\node[below] at (A){$0$};
\coordinate (A') at (8.5, 0);
\node[below] at (A'){$\lambda_{k+1}$};
\coordinate (C) at (13,0);
\node[below] at (C){$\lambda_{k}$};
\coordinate (D) at (18,0);
\coordinate (S1) at (10.75,0);
\node[above] at (S1){$\sigma_p-\frac{\delta_p}{2}$};
\coordinate (A1) at (16,0);
\node[below] at (A1){$\lambda_{1}$};
\coordinate (S2) at (3.6,0);
\node[above] at (S2){$-\sigma_p+\frac{\delta_p}{2}$};
\coordinate (G1) at (2.5,1);
\node[above] at (G1){$\Gamma_2$};
\coordinate (G2) at (14,1);
\node[above] at (G2){$\Gamma_1$};
\coordinate (B) at (10,0);

\coordinate (B') at (10,1);
\coordinate (F) at (12,1);

\coordinate (G) at (12,-1);

\draw[-] (0,0) -- (19,0);

\draw[thick,blue] (A) -- (A'); 
\draw[thick,brown] (A')  -- (C); 

\filldraw (A1) circle (2pt);
\filldraw (A) circle (2pt);
\filldraw (A') circle (2pt);
\filldraw (C) circle (2pt);
\filldraw (B) circle (2pt);

\draw[very thick,red,->] (10,1) -- (12,1); 
\draw[very thick,red] (12,1) -- (18,1); 
\draw[very thick,red,->] (18,1) -- (D); 
\draw[very thick,red] (D) -- (18,-1); 
\draw[very thick,red,->] (18,-1) -- (12,-1); 
\draw[very thick,red] (12,-1) -- (10,-1); 
\draw[very thick,red,->] (10,-1) -- (B'); 
\draw[very thick,red] (B') -- (10,1); 

\coordinate (A'1) at (5, 0);
\node[below] at (A'1){$\lambda_{n-(p-k)}$};
\coordinate (C'1) at (3,0);
\node[below] at (C'1){$\lambda_{n-(p-k)+1}$};

\coordinate (D'1) at (0.5,0);
\coordinate (A11) at (1,0);
\node[below] at (A11){$\lambda_{n}$};

\coordinate (B'1) at (4,0);
\coordinate (B'11) at (4,1);

\filldraw (B'1) circle (2pt);
\draw[thick,blue] (A) -- (A'); 
\draw[thick,brown] (A')  -- (C); 
\draw[thick,blue] (A) -- (A'1);
\draw[thick,brown] (A'1) -- (C'1);
\filldraw (A'1) circle (2pt);
\filldraw (C'1) circle (2pt);
\filldraw (A11) circle (2pt);
\draw[very thick,red,->] (0.5,1) -- (4,1); 
\draw[very thick,red,->] (4,1) -- (4,-1); 
\draw[very thick,red,->] (4,-1) -- (0.5,-1); 
\draw[very thick,red,->] (0.5,-1) -- (0.5,1); 
\end{tikzpicture}}

\vskip2mm

\subsection{Proof of Theorem \ref{cor: rec} using  Theorem \ref{theo: main} } \label{subsec: proofrec}

 We symmetrize $A$ and $\tilde{A}$ as follows. Set  $\mathcal{A}:= \begin{pmatrix}
0 & A \\
A^\top & 0
\end{pmatrix}$ and $ \mathcal{E}:= \begin{pmatrix}
0 & E \\
E^\top & 0
\end{pmatrix}.$
Given the singular decomposition of $A = U \Sigma V^\top$, it is easy to see  that $\mathcal{A}$ admits the following spectral decomposition
$$\mathcal{A}= \begin{pmatrix}
    \frac{U}{\sqrt{2}} & \frac{U}{\sqrt{2}} \\
    \frac{V}{\sqrt{2}} & -\frac{V}{\sqrt{2}}
\end{pmatrix} \begin{pmatrix}
    \Sigma & 0 \\
    0 & -\Sigma
\end{pmatrix} \begin{pmatrix}
    \frac{U^\top}{\sqrt{2}} & \frac{V^\top}{\sqrt{2}} \\
    \frac{U^\top}{\sqrt{2}} & - \frac{V^\top}{\sqrt{2}}
\end{pmatrix}.$$
Therefore, the rank-$2p$ approximation of $\mathcal{A}$ is
$$\mathcal{A}_{2p}:=\begin{pmatrix}
    \frac{U}{\sqrt{2}} & \frac{U}{\sqrt{2}} \\
    \frac{V}{\sqrt{2}} & -\frac{V}{\sqrt{2}}
\end{pmatrix} \begin{pmatrix}
    \Sigma_p & 0 \\
    0 & -\Sigma_p
\end{pmatrix} \begin{pmatrix}
    \frac{U^\top}{\sqrt{2}} & \frac{V^\top}{\sqrt{2}} \\
    \frac{U^\top}{\sqrt{2}} & - \frac{V^\top}{\sqrt{2}}
\end{pmatrix}, $$
which equals 
$$\mathcal{A}_{2p} = \begin{pmatrix}
     0 & U \Sigma_p V^\top \\
    V \Sigma_p U^\top & 0
\end{pmatrix}.$$
Similarly, we have
$$\tilde{\mathcal{A}}_{2p} = \begin{pmatrix}
     0 & \tilde U \tilde\Sigma_p \tilde V^\top \\
    \tilde V \tilde \Sigma_p \tilde U^\top & 0
\end{pmatrix}.$$
Hence,  
$$\tilde{\mathcal{A}}_{2p} - \mathcal{A}_{2p} = \begin{pmatrix}
     0 & \tilde{U} \tilde{\Sigma}_p \tilde{V}^\top- U \Sigma_p V^\top \\
   \tilde{V} \tilde{\Sigma}_p \tilde{U}^\top - V \Sigma_p U^\top & 0
\end{pmatrix}.$$
It yields
 $$ \Norm{\tilde{\mathcal{A}}_{2p} - \mathcal{A}_{2p}}=\Norm{ \tilde{U} \tilde{\Sigma}_p \tilde{V}^\top- U \Sigma_p V^\top} = \big\|\tilde{A}_{p} - A_{p}\| .$$

 Applying Theorem \ref{theo: main} to the pair of symmetric matrices $(\tilde{\mathcal{A}}, \mathcal{A})$, we obtain Theorem \ref{cor: rec}. 
The rest of the paper is devoted to Theorem \ref{theo: mainPositive}. 
 
\section{Proof of Theorem \ref{theo: mainPositive}} \label{sec: proof}

\subsection{Estimating \texorpdfstring{$\Norm{F_1}$}{F1} } \label{subsec: F1}


Recall our strategy that we have split $(zI-A)^{-1}= \sum_{i=1}^{n} \frac{1}{z-\lambda_i}u_i u_i^{\top} $ as 
$$(zI-A)^{-1}= \sum_{i=1}^{n} \frac{1}{z-\lambda_i}u_i u_i^{\top} = P+Q ,$$ where 
 \begin{equation} \label{P} P :=\sum_{i=1}^{r} \frac{u_i u_i^\top}{z-\lambda_i} \,\,\,\text{and}\,\,\, Q :=  \frac{1}{z} \sum_{j>r} { u_j u_j^\top}. 
 \end{equation} 
 
 Set $U_{\perp} :=  \sum_{j>r}  u_j u_j^\top$. By definition,  $U_{\perp} + \sum_{i \le r } u_i u_i^\top = I_n$, the identity matrix. 
Recall that 
$$F_1 = \frac{1}{2 \pi \textbf{i}}\int_{\Gamma}z (P+Q) E(P+Q)dz = \frac{1}{2 \pi \textbf{i}}\int_{\Gamma} z\left( PEP + QEP+PEQ +QEQ  \right) dz.$$

Since $\Gamma$ does not enclose the origin,  $\frac{1}{2 \pi \textbf{i}}\int_{\Gamma} \frac{z}{z^2} dz =\frac{1}{2 \pi \textbf{i}}\int_{\Gamma} \frac{1}{z} dz =0 $. Therefore,
\begin{equation} \label{estimate1} \frac{1}{2 \pi \textbf{i}}\int_{\Gamma}z QEQ  dz =  U_{\perp} E U_{\perp} \times \frac{1}{2 \pi \textbf{i}} \int_{\Gamma}  \frac{1}{z} dz = 0. 
\end{equation}

Next, by Cauchy's integral theorem (Theorem \ref{theo: Cauchy}),
$$\frac{1}{2 \pi \textbf{i}}\int_{\Gamma} zPEQ dz = \frac{1}{2 \pi \textbf{i}} \int_{\Gamma} \sum_{i=1}^{r} \frac{z}{z(z-\lambda_i)} u_i u_i^\top E \sum_{j>r} u_j u_j^\top dz = \sum_{i=1}^p u_i u_i^\top E \sum_{j>r} u_j u_j^\top.$$
It follows that  
\begin{equation}\label{estimate2} 
    \begin{split}
\Norm{\frac{1}{2 \pi \textbf{i}}\int_{\Gamma}z PEQ dz } & = \big\|\sum_{i=1}^p u_i u_i^\top E \sum_{j>r} u_j u_j^\top \big\| \\
& = \big\|\bigg(\sum_{i=1}^p u_i u_i^\top \bigg) \times E \times\bigg( \sum_{j>r} u_j u_j^\top \bigg) \big\| \,\,\\
& \leq  \|\sum_{i=1}^p u_i u_i^\top\| \times \|E\| \times \big\|\sum_{j>r} u_j u_j^\top \big\| = \|E\|.
    \end{split}
\end{equation}
 The same estimate holds for  $\Norm{\frac{1}{2\pi \textbf{i}}\int_{\Gamma} QEP dz}.$


Now we consider the last term $\frac{1}{2 \pi \textbf{i}}\int_{\Gamma} zPEP dz$, which can be expressed as 
$$ \frac{1}{2 \pi \textbf{i}}\int_{\Gamma} \sum_{ 1\le i, j \leq r} \frac{z}{(z -\lambda_i)(z-\lambda_j)}  u_j u_j^\top E u_i u_i^\top =  \frac{1}{2 \pi \textbf{i}} \int_{\Gamma} \sum_{ 1\le i, j \leq r} \frac{z}{(z -\lambda_i)(z-\lambda_j)}  u_j (u_j^\top E u_i ) u_i^\top .$$
By Lemma \ref{contourZ}, if $1 \leq i,j \leq p$, then $\frac{1}{2\pi \textbf{i}}\int_\Gamma \frac{z}{(z- \lambda_i)(z- \lambda_j)} dz = 1$. If both 
$ i,j>p$, then $\frac{1}{2\pi \textbf{i}}\int_\Gamma \frac{z}{(z- \lambda_i)(z- \lambda_j)} dz = 0$. By Lemma \ref{key contour1}, if $1 \leq i \leq p< j \leq r$, then $|\frac{1}{2 \pi \textbf{i}}\int_{\Gamma } \frac{z}{(z -\lambda_i)(z-\lambda_j)} dz |=|\frac{\lambda_i}{\lambda_i -\lambda_j} | \leq \frac{\lambda_p}{\delta_p}.$ Thus, by the triangle inequality, we obtain 
\begin{equation} \label{PEP0}
    \begin{split}
   &\bigg\|\frac{1}{2 \pi \textbf{i}}\int_{\Gamma}z PEP dz\bigg\| \\
   & \leq \big\|\sum_{1 \leq i,j\leq p}   u_j (u_j^\top E u_i ) u_i^\top \big\| + \frac{\lambda_p}{\delta_p} \bigg( \bigg\| \sum_{1 \leq i \leq p < j \leq r} u_j (u_j^\top E u_i ) u_i^\top \bigg\| +\bigg\|\sum_{1 \leq i \leq p < j \leq r} u_i (u_i^\top E u_j ) u_j^\top\bigg\| \bigg).        \end{split}
\end{equation}
Furthermore, using the definition of spectral norm, we have 
\begin{equation} \label{ujEui}
    \begin{split}
    \bigg\| \sum_{1 \leq i \leq p < j \leq r} u_j (u_j^\top E u_i ) u_i^\top\bigg\| & = \max_{\|\textbf{v}\|=\|\textbf{w}\|=1} \textbf{v}^\top   \sum_{1 \leq i \leq p < j \leq r} u_j (u_j^\top E u_i ) u_i^\top \ \textbf{w} \\
    & \leq \max_{\|\textbf{v}\|=\|\textbf{w}\|=1} \bigg|\textbf{v}^\top   \sum_{1 \leq i \leq p < j \leq r} u_j (u_j^\top E u_i ) u_i^\top  \textbf{w}\bigg| \\
    & \leq \max_{\|\textbf{v}\|=\|\textbf{w}\|=1} \sum_{1 \leq i \leq p < j \leq r} \norm{\textbf{v}^\top   u_j} \times \norm{u_j^\top E u_i} \times |u_i^\top \textbf{w}|,
    \end{split}
\end{equation}
where the last inequality follows from the triangle inequality. By  definition, $|u_j^\top E u_i| \leq x$ for all $1 \leq i,j \leq r$, the RHS is at most 
\begin{equation*}
    \begin{split}
\max_{\|\textbf{v}\|=\|\textbf{w}\|=1} \sum_{1 \leq i \leq p < j \leq r} \norm{\textbf{v}^\top   u_j} \times x\times |u_i^\top \textbf{w}| & =   \max_{\|\textbf{v}\|=\|\textbf{w}\|=1} x \bigg(\sum_{p < j \leq r} \norm{\textbf{v}^\top   u_j} \bigg) \bigg( \sum_{1 \leq i \leq p} |u_i^\top \textbf{w}| \bigg). 
    \end{split}
\end{equation*}
Since $u_1,u_2, \dots, u_n $ are orthonormal vectors, using the Cauchy-Schwarz inequality, we have 
$$\sum_{p < j \leq r} \norm{\textbf{v}^\top   u_j}  \leq \sqrt{r-p} \,\,\,\text{and}\,\,\,\sum_{1 \leq i \leq p} |u_i^\top \textbf{w}| \leq \sqrt{p}.$$
Thus, we obtain 
\begin{equation} \label{PEP1}
   \bigg\|\sum_{1 \leq i \leq p < j \leq r} u_j (u_j^\top E u_i ) u_i^\top\bigg\| \leq \sqrt{pr} \, x. 
\end{equation}
By a similar argument, we have
\begin{equation} \label{PEP2}
\begin{split}
    &\bigg\|\sum_{1 \leq i \leq p < j \leq r} u_i (u_i^\top E u_j ) u_j^\top\bigg\| \leq \sqrt{pr} x \,\,\text{and}\,\, \bigg\|\sum_{1 \leq i,j\leq p}   u_j (u_j^\top E u_i ) u_i^\top\bigg\| \leq p  x.
\end{split}    
\end{equation}
Combining \eqref{PEP0}, \eqref{PEP1} and \eqref{PEP2}, we obtain 
\begin{equation}\label{estimate3}
\begin{split}
\Norm{\frac{1}{2 \pi \textbf{i}}\int_{\Gamma}z PEP dz} & \leq  p x + \frac{\lambda_p}{\delta_p} \times 2 \sqrt{pr} x \leq 3 \sqrt{pr} \frac{\lambda_p x}{\delta_p}.
\end{split}
\end{equation}


By \eqref{estimate1}, \eqref{estimate2}, and \eqref{estimate3}, we obtain    
  \begin{equation} \label{F1} 
   \|F_1 \| \le 2 \|E\|   + 3\sqrt{pr} \frac{\lambda_p x}{\delta_p} . 
   \end{equation}

\subsection{Estimating \texorpdfstring{$\|F_s\|$}{Fs}, \texorpdfstring{$s \ge 2$}{s>2}.} Recall the definition 
$$F_s:= \frac{1}{2 \pi \textbf{i}} \int_{\Gamma} z(zI-A)^{-1}[E (zI-A)^{-1}]^s dz.$$
Similar to the previous subsection, we split $(zI-A)^{-1}$ sa $P+Q$, obtaining 
$$F_s:= \frac{1}{2 \pi \textbf{i}} \int_{\Gamma} z(P+Q)[E (P+Q)]^s dz.$$
We  expand $z(P+Q)[ E(P+Q)]^s$ into the sum of $2^{s+1}$ operators, each of which is a product of alternating $Q$-blocks and 
$P$-blocks 
$$z(QEQE \dots QE)(PE PE \dots PE) \dots (PE PE \dots PE) ( QE  QE \dots Q), $$ where we allow the first and last blocks to be empty.

We code each operator like this by the numbers of $Q$'s and the numbers of $P$'s in each block. If there are $(k+1)$  $Q$-blocks and $k$ $P$-blocks, for some integer $k$,  we let 
$\alpha_1, \dots, \alpha_{k+1} $ and $\beta_1, \dots, \beta_k$ be these numbers. These numbers should satisfy the following conditions 
$$ \alpha_1, \alpha_{k+1} \geq 0,$$
$$ \alpha_i, \beta_j \geq 1,\,\, \forall 1<i<k+1,\,\, \forall  1\le j \le k, $$
$$\alpha_1+...+\alpha_{k+1}+\beta_1 +...+\beta_k=s+1.$$

In what follows, we set $\alpha :=(\alpha_1,...,\alpha_{k+1}) \in \mathbb{Z}_{\geq 0}^{k+1} $ and $\beta :=(\beta_1,...,\beta_k) \in \mathbb{Z}_{\geq 0}^{k}$, and use $M(\alpha; \beta)$ to denote the corresponding operator.

\vskip2mm 

\noindent \textit{Example:} Set $s=10, k=2, \alpha_1=3, \alpha_2=2, \alpha_3=1, \beta_1=2, \beta_2=3$; we have 
 $$\int_{\Gamma} M((3,2,1);(2,3)) dz = \int_{\Gamma} z \underset{\substack{\alpha_1=3\\ \text{numbers of $Q$}}}{\underbrace{(QEQEQE)}} \underset{\substack{\beta_1=2\\ \text{numbers of $P$}}} {\underbrace{(PEPE)}}\,\,\, \underset{\substack{\alpha_2=2 \\ \text{numbers of } \, Q}}{\underbrace{(QEQE)}} \underset{\substack{\beta_2=3 \\ \text{numbers of}\, P}}{\underbrace{(PEPEPE)}} \,\,\, \underset{\substack{\alpha_3=1 \\ \text{numbers of}\, Q}}{\underbrace{Q}} dz.$$

Our strategy is to bound
\(
\|\frac{1}{2\pi\mathbf{i}}\int_{\Gamma}M(\alpha;\beta),dz\|
\)
for each $(\alpha;\beta)$ and then sum these bounds using the triangle inequality. We call $(\alpha;\beta)$ \emph{special} if, for some $k\ge2$,
\(
\alpha_1=\alpha_{k+1}=0, \) and
\(
\alpha_2=\cdots=\alpha_k=\beta_1=\cdots=\beta_k=1,
\)
or equivalently,
\[
M(\alpha;\beta)= zPEQEPEQ\cdots PEQEP.
\]
Otherwise, $(\alpha;\beta)$ is called \emph{non-special}.

The following lemma, which bounds
\(
\|\frac{1}{2\pi\mathbf{i}}\int_{\Gamma}M(\alpha;\beta),dz\|,
\)
is the key ingredient of the proof.

\begin{lemma} \label{lemma: Mab} Under the assumption of Theorem \ref{theo: mainPositive}, we have:
\begin{itemize}
\item For any non-special pair $(\alpha;\beta)$, 
\begin{equation*}
  \frac{\Norm{\frac{1}{2 \pi \textbf{i}} \int_{\Gamma} M(\alpha; \beta) dz}}{\lambda_p}  \le  \frac{ 2p }{ 12^s}  \left(  \frac {\|E\|}{\lambda_p} + \frac{rx}{\delta_p} \right).
\end{equation*}
    \item For any special pair ($\alpha;\beta$), 
    \begin{equation*}
        \frac{\Norm{\frac{1}{2 \pi \textbf{i}} \int_{\Gamma} M(\alpha; \beta) dz}}{\lambda_p}  \le   4r \Big(  \frac{ ry +r^2 x^2}{ \lambda_p \delta_p }  \Big) \Big( \frac{2\sqrt{r}\|E\|}{\sqrt{\lambda_p \delta_p}} \Big)^{s-2} +  p \left( \frac{2 \| E\|}{ \lambda_p }  \right)^s.
    \end{equation*}
    \end{itemize}

    \end{lemma}

We prove Lemma \ref{lemma: Mab} in the next two sections.
Assuming Lemma \ref{lemma: Mab}, we can complete the proof of Theorem \ref{theo: mainPositive}. 
For each  $F_s$, there are at most $2^{s+1}$ pairs $(\alpha, \beta)$, in which there is at most one special pair $(\alpha;\beta)$.  Thus, Lemma 
\ref{lemma: Mab} implies that 
 \begin{equation*}
     \begin{split}
  \frac{\| F_s \|}{\lambda_p} & \leq \sum_{(\alpha;\beta)} \frac{\Norm{\frac{1}{2 \pi \textbf{i}} \int_{\Gamma} M(\alpha; \beta) dz}}{\lambda_p} \\
  & \leq  \frac{ 4p }{ 6^s}  \Big(  \frac {\|E\|}{\lambda_p} + \frac{rx}{\delta_p} \Big)   +  4r \Big(  \frac{ ry +r^2 x^2}{ \lambda_p \delta_p }  \Big) \Big( \frac{2\sqrt{r}\|E\|}{\sqrt{\lambda_p \delta_p}} \Big)^{s-2} +   p \left( \frac{2 \| E\|}{ \lambda_p }  \right) ^s := G_s.
     \end{split}
 \end{equation*}
 Moreover, the assumption of Theorem \ref{theo: mainPositive} that  $ \frac{p\|E\|}{\lambda_p}, \frac{\sqrt{r} \|E\|}{\sqrt{\lambda_p \delta_p}}  \leq \frac{1}{24} $ implies  $\frac{2\|E\|}{\lambda_p}, \frac{2\sqrt{r}\|E\|}{\sqrt{\lambda_p \delta_p}} \leq \frac{1}{12}$, yielding that $G_{s+1} \le \frac{1}{6} G_s$, for all $s \ge 2$. Thus,  
 \begin{equation}  \label{F2} \frac{\sum_{s\ge 2} \|F_s\|}{\lambda_p} \leq \sum_{s \geq 2} G_s    
 \le 2 G_2 \leq p\Big(  \frac {\|E\|}{\lambda_p} + \frac{rx}{\delta_p} \Big)  + 8 r  \frac{ ry +r^2 x^2}{ \lambda_p \delta_p } + 8 p  \Big( \frac{ \| E \| }{ \lambda_p } \Big)^2 .
 \end{equation}


 Combining \eqref{bound00}, \eqref{F1}, \eqref{F2}, we obtain 
 \begin{equation} \label{F3}  
 \frac{\| \tilde{A}_p- A_p\|}{\lambda_p} \le 2  \frac{\|E\|}{\lambda_p}   + 3\sqrt{pr} \frac{ x}{\delta_p} + p\Big(  \frac {\|E\|}{\lambda_p} + \frac{rx}{\delta_p} \Big)  + 8r  \frac{ ry +r^2 x^2}{ \lambda_p \delta_p } + 8 p  \Big( \frac{ \| E \| }{ \lambda_p } \Big)^2 . 
 \end{equation} 
 
 \noindent Once again, by the assumption \eqref{ass: sym0} of Theorem \ref{theo: mainPositive}, $\frac{ p\| E \| }{\lambda_p }, \frac{rx}{\delta_p} \le \frac{1}{24} $, so we can replace $8 p  \left(\frac{ \| E \| }{ \lambda_p } \right)^2 $ with $ \frac{ \| E \| }{\lambda_p }$, and  $8\frac{r^3 x^2}{ \lambda_p \delta_p }$ with $ \frac{r x}{ \delta_p }$. 
  Thus, \eqref{F3} simplifies to  
  \begin{equation} \label{F4} 
 \frac{\| \tilde{A}_p- A_p\|}{\lambda_p}  
 \le \frac{(p+3)\| E \| }{ \lambda_p} + \frac{ (p+4) r x } {\delta_p} + 8 \frac{ r^2 y }{ \lambda_p \delta_p }.
 \end{equation} 
 Since $p \geq 1$ and hence $p+4 \leq 5p \leq  8p$, \eqref{F4} simply yields 
 $$\big\|\tilde{A}_p -A_p\big\| \leq 8 \lambda_p \left(\frac{p\|E\|}{\lambda_p} + \frac{pr x}{\delta_p} + \frac{r^2 y}{\lambda_p \delta_p} \right),\,\,\text{proving Theorem \ref{theo: mainPositive}.} $$

 \begin{remark}
     The proof of Theorem \ref{theo: main} follows a similar strategy. It is sufficient to replace Lemma \ref{lemma: Mab} by the following analogue, where 
     $\delta_S$ replaces $\delta_p$ (see \eqref{gapdef} for the definition of $\delta_S$).
 \end{remark}
 \begin{lemma} \label{lemma: MabPi}
   Under the assumption of Theorem \ref{theo: main}, we have
   \begin{itemize}
       \item For any non-special pair ($\alpha;\beta$), 
       \[
       \frac{\Norm{\frac{1}{2 \pi \textbf{i}} \int_{\Gamma} M(\alpha; \beta) dz}}{\sigma_p}  \le \frac{ 2p }{ 12^s}  \left(  \frac {\|E\|}{\sigma_p} + \frac{rx}{\delta_S} \right).
       \]
       \item For any special pair $(\alpha;\beta)$,
       \[
       \frac{\Norm{\frac{1}{2 \pi \textbf{i}} \int_{\Gamma} M(\alpha; \beta) dz}}{\sigma_p}  \le 4r \Big(  \frac{ ry +r^2 x^2}{ \sigma_p \delta_S }  \Big) \Big( \frac{2\sqrt{r}\|E\|}{\sqrt{\sigma_p \delta_S}} \Big)^{s-2} +r \left( \frac{2 \| E\|}{ \sigma_p }  \right) ^s.
       \]
   \end{itemize}
 \end{lemma}

\section{Proof of Lemma \ref{lemma: Mab}: Some preparations} \label{proof: lemma}

In this section, we first derive a few simple contour estimates using Cauchy's integral theorem. Next, we establish our key technical tool, Lemma \ref{key contour1}, using these estimates and an inductive argument. The full proof of 
Lemma \ref{lemma: Mab} is given in the next section. 

\subsection{Preparation: some simple contour estimates }

First, we need the following well-known fact. 

\begin{lemma} \label{contour0} Let $\Gamma$ be a simple closed contour.  If $a$ lies inside the contour, then 
$\frac{1}{2 \pi \textbf{i}}\int_{\Gamma} \frac{1}{ (z-a)^l} dz= 0
\,\, \text{ for all}\,\, l \ge 2.$ If $a$ lies outside the contour, then 
$\frac{1}{2 \pi \textbf{i}}\int_{\Gamma} \frac{1}{ (z-a)^l} dz= 0,\,\, \text{for all}\,\,l \ge 1.$ 
\end{lemma}

\begin{lemma} \label{integral} Let $\Gamma$ be a simple closed contour and let $a, b$ be two points  not on $\Gamma$. Then 
$$ \norm{\frac{1}{2 \pi \textbf{i}}\int_{\Gamma} \frac{1}{(z-a) (z-b) } dz }   =\frac{1 }{  |b-a| } $$ if exactly one of $a$ and $b$ lies inside $\Gamma$. In all other cases, the integral equals zero.  \end{lemma} 
To prove Lemma \ref{integral}, notice that 
\begin{equation} \int_{\Gamma} \frac{1}{(z-a) (z-b)} dz  = \frac{1}{b-a} \Big( \int_{\Gamma} \frac{1}{ z-a } dz - \int_{\Gamma} \frac{1}{ z-b} dz \Big).  \end{equation} 
The claim follows from Lemma \ref{contour0}. Using Lemma \ref{contour0} and Lemma \ref{integral}, we derive two more lemmas. 
\begin{lemma} \label{contour} Let $\Gamma$ be a simple closed contour and $a_1, \dots, a_m$ be points inside $\Gamma$. For any collection of nonnegative integers $(t_1, t_2, \dots, t_m)$, if $\sum_{i=1}^m t_i \ge 2$, then 
$$ \frac{1}{2 \pi \textbf{i}}\int_{\Gamma} \frac{1}{ \prod_{i=1}^m (z-a_i)^{t_i} } dz =0 .$$ 
\end{lemma}

\begin{proof}[Proof of Lemma \ref{contour}]
  If $s:= \sum_{i=1}^m t_i =2$, then the integrand is either $\frac{1}{2 \pi \textbf{i}}\int_{\Gamma} \frac{1}{ (z-a)^2 }dz$ or  $\frac{1}{2 \pi \textbf{i}}\int_{\Gamma} \frac{1}{ (z-a)(z-b)}dz  $ where both $a, b$ are inside the contour. In this case, the claim follows from the previous two lemmas. 
Now we use induction on $s$. Consider $s \ge 3$. If only one of the $t_i $ is positive, then the claim follows again from Lemma \ref{contour0}. If at least two exponents, say $t_i $ and $t_j$,  are positive, then write 
$$ \frac{1} { (z-a_i) (z-a_j)}  = \frac{1}{ a_j -a_i } \Big( \frac{1}{ z-a_i} -\frac{1}{ z-a_j} \Big) . $$
This identity allows us to split the integrand into the sum of two sub-integrands, and in each of these, the sum of exponents equals $s-1$. Applying the induction hypothesis, we conclude the proof.   
\end{proof}

\begin{lemma} \label{contourZ} Let $\Gamma$ be a simple closed contour and $a_1, \dots, a_m$ be points inside $\Gamma$. For any collection of nonnegative integers $(t_1, t_2, \dots, t_m)$, 
\begin{itemize}
    \item if $\sum_{i=1}^m t_i > 2$, then
$  \frac{1}{2 \pi \textbf{i}} \int_{\Gamma} \frac{z}{ \prod_{i=1}^m (z-a_i)^{t_i} } dz = 0 ;$ 
\item if $\sum_{i=1}^m t_i = 2$, then 
$ \frac{1}{2 \pi \textbf{i}}\int_{\Gamma} \frac{z}{ \prod_{i=1}^m (z-a_i)^{t_i} } dz = 1 .$ 
\end{itemize} 

\end{lemma} 
\begin{proof}[Proof of Lemma \ref{contourZ}] Without  loss of generality, we assume that $t_1 > 0$.
    Writing  $z$ as $(z -a_1)+ a_1$, we have 
    \begin{equation*}
        \begin{split}
          \frac{1}{2 \pi \textbf{i}}\int_{\Gamma} \frac{z}{ \prod_{i=1}^m (z-a_i)^{t_i} } dz & = \frac{1}{2 \pi \textbf{i}}\int_{\Gamma} \frac{z-a_1}{ \prod_{i=1}^m (z-a_i)^{t_i} } dz  + \frac{1}{2 \pi \textbf{i}}\int_{\Gamma} \frac{a_1}{ \prod_{i=1}^m (z-a_i)^{t_i} } dz \\
          & = \frac{1}{2 \pi \textbf{i}}\int_{\Gamma} \frac{1}{(z-a_1)^{t_1-1} \prod_{i >1} (z-a_i)^{t_i} } dz  + \frac{a_1}{2 \pi \textbf{i}}\int_{\Gamma} \frac{1}{ \prod_{i=1}^m (z-a_i)^{t_i} } dz.
          \end{split}
          \end{equation*}

By Lemma \ref{contour}, the second term on the RHS is zero, so the LHS equals 
$$ \frac{1}{2 \pi \textbf{i}}\int_{\Gamma} \frac{1}{(z-a_1)^{t_1-1} \prod_{i >1} (z-a_i)^{t_i} } dz. $$

If $\sum_{i=1}^m t_i > 2$ then $t_1 -1+\sum_{i>1} t_i \geq 2$. In this case, again by 
Lemma \ref{contour}, 
$$ \frac{1}{2 \pi \textbf{i}}\int_{\Gamma} \frac{1}{(z-a_1)^{t_1-1} \prod_{i >1} (z-a_i)^{t_i} } dz = 0, $$ which proves the first claim of the lemma. 
 In the case when  $\sum_{i=1}^m t_i =2$,  $(t_1-1)+ \sum_{i>1} t_i=1$. Since all $a_1,a_2,..., a_m$ are inside the contour, 
 $$\frac{1}{2 \pi \textbf{i}}\int_{\Gamma} \frac{1}{(z-a_1)^{t_1-1} \prod_{i >1} (z-a_i)^{t_i} } dz=1, $$  by Cauchy's integral theorem. We complete the proof.    
\end{proof}

\subsection {Our key technical tool }
We first present a generalization of Lemma \ref{contour}. Using this generalization, we then state and prove Lemma \ref{key contour1}, which will be our key technical tool for the proof of Lemma \ref{lemma: Mab}. 



\begin{lemma} \label{key contour0}  Let $\Gamma$ be a simple closed contour and $a_1, \dots, a_l, a_{l+1} , \dots, a_m$, 
$1\le l < m $ be distinct points in the complex plane. Assume that the points in $ X= \{ a_1, \dots, a_l \}$ lie
inside $\Gamma$ and the points in  $Y = \{ a_{l+1} , \dots, a_m \}$ lie outside.   Let $t_i$ be non-negative  integers and $s :=\sum_{i=1}^m t_i$. Furthermore, let $\delta:= \min_{i  \in X, j \in Y } |a_i-a_j| $
and $\lambda = \min_{i \in X} | a_i -a_m |$, and $t:= t_m $. 
Then 
\begin{equation} \label{contourbound0}
\norm{\frac{1}{2 \pi \textbf{i}}\int_ {\Gamma} \frac{1}{ \prod _{i=1}^m (z -a_i)^ {t_i } } dz  } \le  \frac{2^{s-1} }{ \lambda^t \delta ^{s-t-1}} .  
\end{equation}
    \end{lemma}


\begin{remark} The choice of $a_m$ is arbitrary. We can choose any element of $Y$. \end{remark}

\begin{proof}[Proof of Lemma \ref{key contour0}] 
Notice that by definition $\lambda \ge \delta$.     We first prove a weaker bound  that $ \norm{\frac{1}{2 \pi \textbf{i}}\int_ {\Gamma} \frac{1}{ \prod _{i=1}^m (z -a_i)^ {t_i } } dz  } \le  \frac{2^{s-1}}{\delta^{s-1}}$.

We use induction on $s$. For the base case $s=1$, by Cauchy's integral theorem, we have 
$\frac{1}{2 \pi \textbf{i}}\int_ {\Gamma} \frac{1}{ \prod _{i=1} ^m (z -a_i)^ {t_i } } dz $ equals $1$ if $\sum_{a_j \in Y}t_j =0$ and $0$ if $\sum_{a_i \in X} t_i=0$. Thus, the claim holds. 

Now consider $s \ge 2$, and assume that the bound holds for $s-1$.  Notice that if either $\sum_{a_i \in X} t_i = 0$ or $\sum_{a_j \in Y} t_j = 0$, then, by Lemma \ref{contour}, the integral is zero and the bound holds trivially. In the other case, we can  choose 
$a_i \in X$ and $a_j \in Y$ such that both  $t_i, t_j \ge 1$. Rewrite 
$$ \frac{1}{ (z-a_i) (z-a_j) } = \frac{1} {a_j -a_i} \Big( \frac{1}{z-a_i} - \frac{1}{z-a_j } \Big) . $$
By this way, the function $\frac{1}{ \prod _{i=1} ^m (z -a_i)^ {t_i } } dz $ under the integration splits into two parts, 
the first with exponents $$t_1, \dots, t_{i-1} , t_i-1, t_{i+1}, \dots, t_m , $$ and the second with exponents 
$$t_1, \dots, t_{j-1}  , t_j-1, t_{j+1} , \dots, t_m . $$

In both parts,  the sum of the exponents is $s-1$. Now applying the induction hypothesis for both parts 
and  taking into account the fact that the extra factor $|\frac{1}{ a_i -a_j}| \le \frac{1}{\delta} $, we complete the induction.

Now we  prove the stronger bound that $\norm{\frac{1}{2 \pi \textbf{i}}\int_ {\Gamma} \frac{1}{ \prod _{i=1}^m (z -a_i)^ {t_i } } dz  } \le \frac{2^{s-1}}{\lambda^{t} \delta ^{s-t-1}}. $ Again, we use induction on $s$. The case $s=1$ is the same as above. For $s \ge 2$, we consider two cases. If $t=0$, then we are done by the previous consideration, as there is no $\lambda$ in the bound. If $t \ge 1$, then we split 
$$ \frac{1}{ (z-a_i) (z-a_m) } = \frac{1} {a_m -a_i} \left(\frac{1}{z-a_i} - \frac{1}{z-a_m} \right) . $$

Applying the induction hypothesis in each part, and taking into account that $|\frac{1} {a_i -a_m}|  \le \lambda^{-1} $, we obtain the bound 
$$ \lambda^{-1} \left(\frac{2^{s-2}} 
{  \lambda ^{t-1} \delta ^{[(s-1) -(t-1)]-1  } }  +  \frac{  2^{s-2}}  {  \lambda ^{t} \delta ^{[(s-1 ) -t ]-1} } \right) \le  \frac{ 2^{s-1} } { \lambda ^{t} \delta ^{s-1-t} }  , $$ completing  the proof.
\end{proof}

Now we state and prove our key technical lemma.

\begin{lemma} \label{key contour1}  Let $\Gamma$ be a simple closed contour and $a_1, \dots, a_l, a_{l+1} , \dots, a_m$, 
$1\le l < m$ be distinct points in the complex plane such that the points in $ X= \{ a_1, \dots, a_l \}$ lie 
inside $\Gamma$ and the points in $Y = \{ a_{l+1}, \dots, a_m \}$ lie outside.  Let $t_i$ be positive integers for $1 \leq i \leq m-1$ and let $t_m:=t$ be a nonnegative integer. Set $s :=\sum_{i=1}^m t_i$. Furthermore, define $\delta:= \min_{i  \in X, j \in Y } |a_i-a_j| $
and $\lambda := \min_{i \in X} | a_i -a_m |$. Then,  
\begin{equation} \label{key contourbound1}
\begin{split}
 &  \norm{\frac{1}{2 \pi \textbf{i}} \int_ {\Gamma} \frac{z}{ \prod_{i=1} ^m (z -a_i)^ {t_i } } dz  } = \norm{\frac{a_1}{(a_1 -a_m)^t}} \leq \frac{1}{\lambda^{t-1}} + \frac{|a_{m}|}{\lambda^t}\,\,\,\text{if $s = t+1$},\\
 &  \norm{\frac{1}{2 \pi \textbf{i}} \int_ {\Gamma} \frac{z}{ \prod _{i=1} ^m (z -a_i)^ {t_i } } dz  } \le 2^s \left(\frac{1}{ \lambda^t} \frac{1}{ \delta ^{s-t-2}} + \frac{|a_{l+1}|}{\lambda^t \delta^{s-t-1}}\right)\,\,\text{if $s \geq t+2$}. 
 \end{split}
\end{equation}

\end{lemma} 

\begin{proof}[Proof of Lemma \ref{key contour1}] 

Let us first consider the easier case when  $s=t+1$. In this case, as $t_i \ge 1$ for all $i \le m-1$, $\sum_{a_i \in X} t_i\geq |X|$ and $\sum_{a_j \in Y} t_j \geq t_m \geq t$, we must have $X=\{a_1\}$ with $t_1=1$, and  $ Y=\{a_m\}$ with $t_m=t$. Therefore, 
$$   \prod_{i=1}^m (z -a_i)^ {t_i }  =   (z- a_{1}) (z -a_m)^ {t }  . $$
Writing  $z= (z-a_1) +a_1$, we split the integrand into a sum of two terms, 
$$\frac{1}{2 \pi \textbf{i}} \int_ {\Gamma} \frac{z}{ \prod_{i=1} ^m (z -a_i)^ {t_i } } dz  = \frac{1}{2 \pi \textbf{i}} \int_ {\Gamma} \frac{1}{ (z -a_m)^ {t } } dz
+ \frac{1}{2 \pi \textbf{i}}\int_ {\Gamma} \frac{a_1 }{ (z- a_{1}) (z -a_m)^ {t } } dz. $$
By the triangle inequality 
$$ \norm{ \frac{1}{2 \pi \textbf{i}}\int_ {\Gamma} \frac{z}{ \prod_{i=1} ^m (z -a_i)^ {t_i } } dz  } \le \norm{\frac{1}{2 \pi \textbf{i}} \int_ {\Gamma} \frac{1}{ (z -a_m)^ {t } } dz} 
+ \norm{\frac{1}{2 \pi \textbf{i}} \int_ {\Gamma} \frac{a_1 }{ (z- a_{1}) (z -a_m)^ {t } } dz}. $$
By Cauchy's integral theorem (Theorem~\ref{theo: Cauchy}), $ \frac{1}{2 \pi \textbf{i}}\int_ {\Gamma} \frac{1}{ (z -a_m)^ {t } } dz=0$. Furthermore, 
$$ \norm{\frac{1}{2 \pi \textbf{i}} \int_ {\Gamma} \frac{a_1 }{ (z- a_{1}) (z -a_m)^ {t } } dz} = 
\norm{\frac{a_1}{(a_1 - a_m)^t}}. $$
Notice that in  this case  $\lambda =\delta =|a_1 - a_m|$. Splitting  $a_1= a_m + (a_1-a_m) $, we obtain the desired bound.

\vskip2mm 

For the more difficult case when  $s \geq t+2$, we divide our proof into two sub-cases:  $l+1 < m$ (equivalently $|Y| \geq 2$ and $\sum_{a_j \in Y}t_j \geq t+1$); and  $l+1 =m$ (equivalently $Y =\{a_m\}$). 

\noindent For the first sub-case when $l+1 < m$, we split the term  $z$ on the numerator as $(z-a_{l+1}) + a_{l+1}$. Using  the triangle inequality, we obtain 
\begin{equation*}
   \begin{split}
       & \norm{ \frac{1}{2 \pi \textbf{i}}\int_ {\Gamma} \frac{z}{ \prod _{i=1} ^m (z -a_i)^ {t_i } } dz  } \leq \\
       & \norm{\frac{1}{2 \pi \textbf{i}} \int_ {\Gamma} \frac{1}{ (z-a_{l+1})^{t_{l+1}-1}\prod_{i \neq l+1} (z -a_i)^ {t_i } } dz  } + |a_{l+1}| \norm{ \frac{1}{2 \pi \textbf{i}}\int_ {\Gamma} \frac{1}{ \prod _{i=1} ^m (z -a_i)^ {t_i } } dz  }.
   \end{split} 
\end{equation*}
By assumption, $t_{l+1} \geq 1$ and hence $t_{l+1} -1 \geq 0$. We  can now apply Lemma \ref{key contour0} (with the same sets $X$ and $Y$)  to obtain 
\begin{equation} \label{eq 1keycontour1}
    \norm{\frac{1}{2 \pi \textbf{i}} \int_ {\Gamma} \frac{1}{ (z-a_{l+1})^{t_{l+1}-1}\prod_{i \neq l+1} (z -a_i)^ {t_i } } dz  } \leq  \frac{2^{s-2}}{\lambda^t \delta^{(s-1)-t-1}} =  \frac{2^{s-2} }{\lambda^t \delta^{s-t-2}}.
\end{equation}
Similarly, by Lemma \ref{key contour0}, we also have 
\begin{equation} \label{eq 2keycontour1}
   |a_{l+1}| \norm{ \frac{1}{2 \pi \textbf{i}} \int_ {\Gamma} \frac{1}{ \prod _{i=1} ^m (z -a_i)^ {t_i } } dz  } \leq |a_{l+1}| \times \frac{2^{s-1}}{\lambda^{t}\delta^{s-t-1}}. 
\end{equation}
Combining \eqref{eq 1keycontour1} and \eqref{eq 2keycontour1}, we  obtain 
$$\norm{ \frac{1}{2 \pi \textbf{i}} \int_ {\Gamma} \frac{z}{ \prod _{i=1} ^m (z -a_i)^ {t_i } } dz  } \leq 2^{s-1} \left(\frac{1}{ \lambda^t} \frac{1}{ \delta ^{s-t-2}} + \frac{|a_{l+1}|}{\lambda^t \delta^{s-t-1}} \right), $$
which proves the first sub-case.

In the second sub-case that $l+1=m$,  $Y=\{a_m\}$. In  this setting, $\delta =\lambda =\min_{i\in X}|a_i-a_m|.$ Without  loss of generality, assume that $\lambda=|a_1-a_m|$. It follows that  $|a_1| \leq |a_m| +\lambda$ by the triangle inequality. We split the term $z$ in the numerator as $(z -a_1)+ a_1$. Using the triangle inequality, we obtain 
\begin{equation*}
    \begin{split}
        & \norm{ \frac{1}{2 \pi \textbf{i}}\int_ {\Gamma} \frac{z}{ \prod _{i=1} ^m (z -a_i)^ {t_i } } dz  } \leq \\
        & \norm{ \frac{1}{2 \pi \textbf{i}}\int_ {\Gamma} \frac{1}{ (z-a_{1})^{t_{1}-1}\prod_{i \neq 1} (z -a_i)^ {t_i } } dz  } + |a_{1}| \norm{ \frac{1}{2 \pi \textbf{i}}\int_ {\Gamma} \frac{1}{ \prod _{i=1} ^m (z -a_i)^ {t_i } } dz  }.
    \end{split}
\end{equation*}
Now we apply Lemma \ref{key contour0} to bound both terms.  Note that $\lambda=\delta$ in this sub-case, we obtain 
$$\norm{\frac{1}{2 \pi \textbf{i}} \int_ {\Gamma} \frac{1}{ (z-a_{1})^{t_{1}-1}\prod_{i \neq 1} (z -a_i)^ {t_i } } dz  } \leq  2^{s-2} \frac{1}{\lambda^{s-2}}, $$
and 
$$|a_{1}| \norm{ \frac{1}{2 \pi \textbf{i}}\int_ {\Gamma} \frac{1}{ \prod _{i=1} ^m (z -a_i)^ {t_i } } dz  } \leq |a_1| \times \frac{ 2^{s-1}}{\lambda^{s-1}} \leq (\lambda +|a_m|) \times \frac{ 2^{s-1}}{\lambda^{s-1}}.$$
It follows that 
$$\norm{ \frac{1}{2 \pi \textbf{i}} \int_ {\Gamma} \frac{z}{ \prod _{i=1} ^m (z -a_i)^ {t_i } } dz  } \leq  2^{s-2} \frac{1}{\lambda^{s-2}}+ (\lambda +|a_m|) \times \frac{2^{s-1}}{\lambda^{s-1}} \leq \frac{2^s}{\lambda^{s-2}}+|a_m| \frac{2^{s}}{\lambda^{s-1}},$$
proving the desired bound. The proof of the lemma is completed. 
\end{proof} 


\section{Proof of Lemma \ref{lemma: Mab}} 
\label{proof:lemma2}

Fix an admissible pair $(\alpha, \beta)$. Recall  the definitions that $U_{\perp }= \sum_{i >r} u_{i} u_i^\top $,  
$P= \sum_{i=1}^r  \frac{1}{z -\lambda_i } u_i u_i^\top $ and $Q= \sum_{j >r} \frac{1}{z} u_j u_j^\top  =\frac{1}{z} U_{\perp} $.  Let $s_1=  \sum_{i=1}^{k+1} \alpha_i$ and  $s_2 = \sum_{j=1}^k \beta_j $. Define
\begin{equation} \label{A_hB_hdefinition}
    \begin{split}
        & B_{h} := \beta_1 + \beta_2 + ... + \beta_{h} \,\,\,\, \text{for each} \,\,\, 1 \leq h \leq k, \\
    &  B_0 = 0.
    \end{split}
\end{equation}

Consider a term $M(\alpha; \beta)$ as defined in the lemma. 
Multiplying out all the terms in $P$ (there are $s_2$ of $P$'s), we obtain 
\begin{equation*}
    \begin{split}
&  M(\alpha; \beta)= \sum_{i_1,i_2,...,i_{s_2} \leq r} \frac{z} {z^{s_1} \prod_{j=1}^{s_2}(z-\lambda_{i_j})} \times \\ 
&\left[ \prod_{h=1}^k \left(U_{\perp} E \right)^{\alpha_h} \times \left(\prod_{l=1}^{\beta_{h}} u_{i_{l+B_{h-1}}}u_{i_{l+B_{h-1}}}^\top E \right)   \right] \times \left[ U_{\perp}\prod_{l=1}^{\alpha_{k+1-1}} (E U_{\perp}) \right].
    \end{split}
\end{equation*}


By the triangle inequality,  
\begin{equation*}
    \begin{split}
&  \Norm{\frac{1}{2 \pi \textbf{i}}\int_{\Gamma} M(\alpha; \beta) dz} \leq \sum_{i_1,i_2,...,i_{s_2} \leq r} \norm{\frac{1}{2 \pi \textbf{i}} \int_{\Gamma} \frac{z} {z^{s_1} \prod_{j=1}^{s_2}(z-\lambda_{i_j})} dz} \times W,\,\,\text{where} 
    \end{split}
\end{equation*}
\begin{equation}\label{long}
   W:=\Norm{\left[ \prod_{h=1}^k \left(U_{\perp} E \right)^{\alpha_h} \times \left(\prod_{l=1}^{\beta_{h}} u_{i_{l+B_{h-1}}}u_{i_{l+B_{h-1}}}^\top E \right)   \right] \times \left[ U_{\perp}\prod_{l=1}^{\alpha_{k+1-1}} (E U_{\perp}) \right]}. 
\end{equation}


For each particular sequence $[i_1,i_2, \dots, i_{s_2}]$, we bound $\norm{\frac{1}{2 \pi \textbf{i}} \int_{\Gamma} \frac{z} {z^{s_1} \prod_{j=1}^{s_2}(z-\lambda_{i_j})} dz}$ using Lemma \ref{key contour1}.  Let $I$ denote the set of distinct indices appearing in $[i_1, i_2, i_3, \dots, i_{s_2}]$. (For example, if $s_2=5
$  and $ [i_1,i_2, i_3, i_4, i_5] =[1,1,1,5,5]$, then $I=\{1,5\}.$) Define
$$X:=\{\lambda_{i}:\,\, i \in I\,\,\text{and}\,\,i \leq p\}, \,\,\text{and}\,\, Y:=\{\lambda_{i}:\,\, i \in I\,\,\text{and}\,\,i > p\} \cup \{0\}.$$
We now apply Lemma \ref{key contour1} to these sets $X$ and $Y$. Notice that 
\begin{itemize}
    \item for each $i \in I$, the exponent of $\lambda_i$ equals its multiplicity in $(\lambda_{i_1}, \lambda_{i_2},\dots,\lambda_{i_{s_2}})$,
    \item the point $0 \in Y$ plays the role of $a_m$ with exponent $t_m = s_1$, 
    \item $\lambda =\min_{x \in X}|x -0| = \min_{i \in I, i\leq p}|\lambda_i| \geq \lambda_p,$
    \item  $\delta= \min_{x \in X, y \in Y}|x-y|=\min_{i,j \in I, i \leq p < j}|\lambda_i -\lambda_j| \geq \delta_p.$
\end{itemize}
Since the total sum of exponents here is $s_1+s_2 =s+1$, Lemma \ref{key contour1} yields 
\begin{equation*}  
 \norm{\frac{1}{2 \pi \textbf{i} }\int_{\Gamma} \frac{z}{z^{s_1} \prod_{j=1}^{s_2}(z-\lambda_{i_j})} dz}  \le 2^{s+1} \left(\frac{1}{\lambda^{s_1} \delta^{(s+1)-s_1-2}}+ \frac{\max_{\lambda_i \in Y}|\lambda_i|}{\lambda^{s_1} \delta^{(s+1)-s_1-1}}  \right). 
\end{equation*} 
 This is at most $ 2^{s+1} \left(\frac{1}{\lambda_p^{s_1} \delta_p^{s_2-2}}+ \frac{\max_{\lambda_i \in Y}|\lambda_i|}{\lambda_p^{s_1} \delta_p^{s_2-1}}  \right)$ as $\lambda \geq \lambda_p, \delta \geq \delta_p$ and $(s+1)-s_1=s_2$.
 
Moreover, as all eigenvalues of $A$ outside $\Gamma$ are at most $\lambda_{p+1}$ and $A$ is positive semi-definite, we have $\max_{\lambda_i \in Y} |\lambda_i| \leq \lambda_{p+1}$.  
Thus,
\begin{equation}\label{contour2} 
  \norm{\frac{1}{2 \pi \textbf{i} }\int_{\Gamma} \frac{z}{z^{s_1} \prod_{j=1}^{s_2}(z-\lambda_{i_j})} dz}  \le  2^{s+1} \left(\frac{1}{\lambda_p^{s_1} \delta_p^{s_2-2}}+ \frac{\lambda_{p+1}}{\lambda_p^{s_1} \delta_p^{s_2-1}}  \right) = \frac{ 2^{s+1}\lambda_{p} } { \lambda_p^{s_1} \delta_p ^{s_2 -1 } } .  
\end{equation}

Next, we bound the matrix norm $W$ (for a fixed sequence $i_1,\dots, i_{s_2} $). To start, we write the long product in \eqref{long} as follows. Let $v :=  U_{\perp}E \cdots U_{\perp} E u_{i_1}$ and $w^\top := u_{i_{s_2}}^\top EU_{\perp} \cdots EU_{\perp}$, and 
$$ \gamma := u_{i_1}^\top E \cdots E u_{i_{\beta_1}} u_{i_{\beta_1}}^\top E  U_{\perp} E  U_{\perp} E \cdots U_{\perp}E \cdots E u_{i_{s_2-\beta_k+1}}u_{i_{s_2-\beta_k+1}}^\top E \cdots E u_{i_{s_2}}. $$

Notice that $v, w$ are vectors and $\gamma$ is a scalar; 
the matrix in question can be written as 
$v \gamma w^\top$, and its norm $W$ is bounded by $\| v \| \times \|w\| \times |\gamma |. $

As $u_i$ are unit vectors and $\| EU _{\perp} \| \le \| E\|$, it follows that $\| v\| \le \| E\| ^{\alpha_1} $ and $\|w\| \le \| E\| ^{\alpha_{k+1} } $. It now remains to bound $\gamma$.  
We start to group the terms in $\gamma$ (going from left to right) as follows 
\begin{equation} \label{block1} \gamma = (u_{i_1}^\top E u_{i_2} ) (u_{i_2} ^\top E u_{i_3} ) \dots   \end{equation} 

The first block which is not of the form $(u_{i_j} ^\top E u_{i_{j+1}})$ is 
\begin{equation} \label{block2} u_{i_{\beta_1}}^\top E  U_{\perp} E  \cdots U_{\perp} E u_{i_{\beta_1 +1}} . \end{equation} 

We will have $k-1$ blocks of this type, alternating with the blocks formed by the product of a few $(u_i^\top E u_j)$'s.  Notice that there are totally  $\sum_{i=1} ^k (\beta_i -1) $ terms of the form 
$u_i^\top E u_j$. Since each term   $| u_i^\top E u_j | \le x$ (by the definition of $x$ in \eqref{def: xyz}), we obtain a contribution of $x^{ (\sum_{i=1}^k \beta_i) - k }$.
 To bound the absolute value of the block in \eqref{block2}, just notice that 
$$  |u_i^\top [EU_{\perp}]^{l} E u_j| \leq \|E\|^{l+1}. $$
The resulting total exponent of $\| E\|$ is  $\sum_{i=2}^{k} (\alpha _i +1) $. Putting everything together, we get 
\begin{equation} \label{Wbound}
    W \leq \| E\| ^{\alpha_1 + \alpha_{k+1} } \| E\| ^{\sum_{i=2}^k (\alpha_i +1) } x^{ \sum_{i=1}^k (\beta_i-1)} = \|E\|^{(\sum_{i=1}^{k+1} \alpha_i)+k-1} x^{(\sum_{i=1}^k \beta_i)-k}.
\end{equation}

Combining \eqref{contour2} and \eqref{Wbound}, for each sequence $[i_1,i_2,\dots, i_{s_2}]$, we have 
$$\norm{\frac{1}{2 \pi \textbf{i}} \int_{\Gamma} \frac{z} {z^{s_1} \prod_{j=1}^{s_2}(z-\lambda_{i_j})} dz} \times W \leq \frac{ 2^{s+1}\lambda_{p} } { \lambda_p^{s_1} \delta_p ^{s_2 -1 } } \times \|E\|^{(\sum_{i=1}^{k+1} \alpha_i)+k-1} x^{(\sum_{i=1}^k \beta_i)-k}.$$
This estimate and \eqref{long} imply
$$ \Norm{\frac{1}{2 \pi \textbf{i}} \int_\Gamma M(\alpha; \beta) dz } \le \sum_{i_1,i_2,...,i_{s_2} \leq r} \frac{ 2^{s+1}\lambda_{p} } { \lambda_p^{s_1} \delta_p ^{s_2 -1 } }  \| E\| ^{\alpha_1 + \alpha_{k+1} } \| E\| ^{\sum_{i=2}^k (\alpha_i +1) } x^{ (\sum_{i=1}^k \beta_i) - k },$$
Since the number of sequences $[i_1,i_2, \dots, i_{s_2}]$ in the definition of $M(\alpha;\beta)$ is at most $p\, r^{s_2-1}$ (at least one $\lambda_i$ must be inside the contour), we further obtain
$$\Norm{\frac{1}{2 \pi \textbf{i}} \int_\Gamma M(\alpha; \beta) dz } \le p \, r^{s_2-1} \frac{ 2^{s+1}\lambda_{p} } { \lambda_p^{s_1} \delta_p ^{s_2 -1 } }  \| E\| ^{\alpha_1 + \alpha_{k+1} } \| E\| ^{\sum_{i=2}^k (\alpha_i +1) } x^{ (\sum_{i=1}^k \beta_i) - k }.$$



Dividing both sides by $\lambda_p$ and rearranging the terms on the RHS, we 
obtain the following nicer inequality 
 \begin{equation} \label{Mab-bound}
 \frac{\|\frac{1}{2 \pi \textbf{i}}\int_\Gamma M(\alpha; \beta) dz \|}{ \lambda_{p}} \leq 2^{s+1}p \left( \frac{\|E\|}{\lambda_p} \right)^{T_1} \left( \frac{ rx}{\delta_p} \right)^{T_2} \left( \frac{ \sqrt r \|E\|}{\sqrt {\lambda_p \delta_p}  } \right)^{T_3},
  \end{equation}
where 
$$T_1= \alpha_1+ \alpha_2+ \dots +\alpha_{k+1} - (k-1),$$
$$T_2= \beta_1+\beta_2+\dots+\beta_k - k,$$
$$T_3= 2(k-1).$$


If $T_1+T_2 \ge 1$, then by the assumption \eqref{ass: sym0} that 
$$\max \left\lbrace  \frac{\|E\|}{\lambda_p}, \frac{rx}{\delta_p} , \frac{ \sqrt r \|E\|}{\sqrt { \lambda_p \delta_p }} \right\rbrace   \le \frac{1}{24} ,$$ we conclude that the  RHS of \eqref{Mab-bound} is at most
\begin{equation} \label{goodcasebound} \max \left\lbrace  \frac {\|E\|}{\lambda_p} ,  \frac{rx}{\delta_p} \right\rbrace \times  \frac{  2^{s+1} }{ 24^{T_1+T_2+T_3 -1} } . \end{equation} 

Since $T_1+T_2 +T_3 -1 = s $, it follows that 
\begin{equation} \label{goodbound1} \frac{\|\frac{1}{2 \pi \textbf{i}}\int_\Gamma M(\alpha; \beta) dz \|}{\lambda_{p}}  \le \frac{ 2 }{ 12^s} p \left(  \frac {\|E\|}{\lambda_p} + \frac{rx}{\delta_p} \right) . \end{equation} 
This proves the first part of Lemma \ref{lemma: Mab}.

\subsection{Handling the exceptional case \texorpdfstring{$T_1+T_2=0$}{T1T2=0}} \label{subsec: special case}
Now we consider the case that $T_1+T_2=0$. This means $T_1 =T_2=0$, which can only occur if  $\alpha_1=\alpha_{k+1}=0$ and $\alpha_2=...=\alpha_k=\beta_1=...=\beta_k=1$ for some $k \geq 2$, or, equivalently,
$$ M(\alpha; \beta)= zPEQEPEQEPEQ \dots PEQEP;\,\text{i.e., $(\alpha;\beta)$ is special.}  $$
In this particular case, note that $s+1 = \sum_{i=1}^{k+1} \alpha_i + \sum_{j=1}^k \beta_j = 2k -1$. It implies $2(k-1)=s.$ 

 Again, write $P =\sum _{i \le r} \frac{1}{ z -\lambda_i } u_i u_i^\top $ and $Q = \frac{1}{z} U_{\perp}$. Multiplying out all the terms, we can express 
$zPEQEPEQ \cdots PEQEP$ as $\Sigma_1+\Sigma_2$, where 
$$ \Sigma_1 :=  \sum_{\substack{i_1,i_2,...,i_k \leq r \\ (i_1,...,i_k) \neq  (j,...,j) \\ \forall 1 \leq j \leq r}}   \frac{z \times u_{i_1} u_{i_k}^\top \prod_{l=1}^{k-1} (u_{i_l}^\top E U_{\perp} E u_{i_{l+1}}) }{z^{k-1}   \prod_{l=1}^k (z -\lambda_{i_l} )  }   = \sum_{\substack{i_1,i_2,...,i_k \leq r \\ (i_1,...,i_k) \neq  (j,...,j) \\ \forall 1 \leq j \leq r}} \frac{u_{i_1} u_{i_k}^\top \prod_{l=1}^{k-1} (u_{i_l}^\top E U_{\perp} E u_{i_{l+1}})}{z^{k-2}   \prod_{l=1}^k (z -\lambda_{i_l} ) }.$$ and 
$$
\Sigma_2 := \sum_{j=1}^{r}  \frac{z \times u_j u_j^\top (u_{j}^\top E U_{\perp} E u_{j}) ^{k-1} }{z^{k-1} (z-\lambda_j)^{k}} = \sum_{j=1}^{r}  \frac{ u_j u_j^\top (u_{j}^\top E U_{\perp} E u_{j}) ^{k-1} }{z^{k-2} (z-\lambda_j)^{k}}.
 $$

 In words, $\Sigma_2$ is the sum of the diagonal terms and $\Sigma_1$
 is the sum of the rest. By the triangle inequality, we have
 \begin{equation} \label{PEQEPsplit}
   \Norm{\frac{1}{2 \pi \textbf{i}}\int_{\Gamma}z PEQEPEQ \cdots PEQEP dz} = \Norm{\frac{1}{2 \pi \textbf{i}} \int_{\Gamma} (\Sigma_1 + \Sigma_2) dz}  \le \Norm{\frac{1}{2 \pi \textbf{i}} \int_{\Gamma} \Sigma_1  dz}+ \Norm{\frac{1}{2 \pi \textbf{i}} \int_{\Gamma} \Sigma_2 dz}. 
 \end{equation}
 Applying the triangle inequality again, we  obtain
\begin{equation}\label{PEQEPSigmasplit}
 \begin{split}
 \Norm{\frac{1}{2 \pi \textbf{i}} \int_{\Gamma} \Sigma_1  dz} & = \bigg\| \frac{1}{2 \pi \textbf{i}} \int_{\Gamma}\sum_{\substack{i_1,i_2,...,i_k \leq r \\ (i_1,...,i_k) \neq  (j,...,j) \\ \forall 1 \leq j \leq r}} \frac{u_{i_1} u_{i_k}^\top \prod_{l=1}^{k-1} (u_{i_l}^\top E U_{\perp} E u_{i_{l+1}})}{z^{k-2}   \prod_{l=1}^k (z -\lambda_{i_l} ) } dz \bigg\| \\
 & \leq \sum_{\substack{i_1,i_2,...,i_k \leq r \\ (i_1,...,i_k) \neq  (j,...,j) \\ \forall 1 \leq j \leq r}} \bigg\|\frac{1}{2 \pi \textbf{i}} \int_{\Gamma} \frac{u_{i_1} u_{i_k}^\top \prod_{l=1}^{k-1} (u_{i_l}^\top E U_{\perp} E u_{i_{l+1}})}{z^{k-2}   \prod_{l=1}^k (z -\lambda_{i_l} ) } dz \bigg\| \leq S_1,
 \end{split}   
\end{equation}
 where
 $$S_1:=    \sum_{\substack{i_1,i_2,...,i_k \leq r \\ (i_1,...,i_k) \neq  (j,...,j) \\ \forall 1 \leq j \leq r}}  \norm{\frac{1}{2 \pi \textbf{i}}\int_{\Gamma} \frac{1}{z^{k-2}   \prod_{l=1}^k (z -\lambda_{i_l} )  } dz } \times \|u_{i_1} u_{i_k}^\top\| \times   \prod_{l=1}^{k-1} | u_{i_l}^\top E U_{\perp} E u_{i_{l+1}}|.$$

To bound $\Norm{\frac{1}{2 \pi \textbf{i}} \int_{\Gamma} \Sigma_2 dz} $, we first notice that 
by Lemma \ref{contour}, 
$ \frac{1}{2\pi\textbf{i}}\int_{\Gamma} \frac{1}{z^{k-2} (z- \lambda_j)^k} dz = 0 \,\,\,\text{for all $j > p$}.$
 Thus, $\frac{1}{2 \pi \textbf{i}} \int_{\Gamma} \Sigma_2 dz$ simplifies to
 $$\frac{1}{2 \pi \textbf{i}} \int_{\Gamma} \sum_{j=1}^p \frac{ u_j u_j^\top (u_{j}^\top E U_{\perp} E u_{j}) ^{k-1} }{z^{k-2} (z-\lambda_j)^{k}} dz.$$ Using the triangle inequality, we have 
 
 \begin{equation} \label{PEPEQS2split}
 \begin{split}
     \Norm{\frac{1}{2 \pi \textbf{i}}  \int_{\Gamma} \Sigma_2 dz} & = \Norm{\frac{1}{2 \pi \textbf{i}} \int_{\Gamma} \sum_{j=1}^p \frac{ u_j u_j^\top (u_{j}^\top E U_{\perp} E u_{j}) ^{k-1} }{z^{k-2} (z-\lambda_j)^{k}} dz }  \leq S_2,   
 \end{split}
   \end{equation}
  where
 $$S_2:= \sum_{j=1}^{p}  \norm{\frac{1}{2 \pi \textbf{i}}\int_{\Gamma} \frac{1}{z^{k-2} (z-\lambda_j)^{k}} dz}  \times \|u_j u_j^\top\| \times  |  u_{j}^\top E U_{\perp} E u_{j}|^{k-1}.$$
It remains to bound $S_1$ and $S_2$.

\vskip2mm 

\noindent\textit{Bounding $S_1$.} We bound the terms in each summand of $S_1$ from left to right. First, for any particular sequence $[i_1,i_2, \dots, i_{k}]$, we can bound $\norm{\frac{1}{2 \pi \textbf{i}}\int_{\Gamma} \frac{1}{z^{k-2}   \prod_l (z -\lambda_{i_l} )  } dz }$ using Lemma \ref{key contour0}.  Let $I$ denote the set of distinct indices appearing in $[i_1, i_2, i_3, \dots, i_{k}]$. Define
$$X:=\{\lambda_{i}:\,\, i \in I\,\,\text{and}\,\,i \leq p\}, \,\,\text{and}\,\, Y:=\{\lambda_{i}:\,\, i \in I\,\,\text{and}\,\,i > p\} \cup \{0\}.$$
We now apply Lemma \ref{key contour0} for these sets $X$ and $Y$. Notice that 
\begin{itemize}
    \item for each $i \in I$, the exponent of $\lambda_i$ equals to its multiplicity in $(\lambda_{i_1}, \lambda_{i_2},\dots,\lambda_{i_{k}})$,
    \item the point $0 \in Y$ plays the role of $a_m$ with exponent $t_m = k-2$, 
    \item $\lambda =\min_{x \in X}|x -0| = \min_{i \in I, i\leq p}|\lambda_i| \geq \lambda_p.$
    \item  $\delta= \min_{x \in X, y \in Y}|x-y|=\min_{i,j \in I, i \leq p < j}|\lambda_i -\lambda_j| \geq \delta_p.$
\end{itemize}
Since the total sum of exponents here is $2k-2= s$, Lemma \ref{key contour0} yields 
\begin{equation} \label{S1integral}  
 \norm{\frac{1}{2 \pi \textbf{i}}\int_{\Gamma} \frac{1}{z^{k-2}   \prod_l (z -\lambda_{i_l} )  } dz } \le \frac{ 2^s}{ \lambda^{k-2} \delta^{k-1} } \leq \frac{2^s}{\lambda_p^{k-2} \delta_p^{k-1}}  =  \frac{ 2^s \lambda_p}{ \lambda_p^{s/2} \delta_p^{ s/2 } }. 
\end{equation} 

Next, we bound  $\prod_{l=1}^{k-1} | u_{i_l}^\top E U_{\perp} E u_{i_{l+1}}|$. For any $1 \leq i,j \leq n$,
 $$| u_{i}^\top E U_{\perp} E u_{j}| \leq \|E\|^2,$$
 since $\|U_{\perp}\|=1$.
For $i\ne j$, we sharpen this bound using the parameters $x$ and $y$ (see \eqref{def: xyz}). 
With $U_{\perp}=I-\sum_{k\le r}u_k u_k^\top$, we write
 $$   u_{i}^\top E U_{\perp} E u_{j} = Eu_i \cdot Eu_j - \sum_{k \le r} (u_i^\top E u_k) ( u_k^\top E u_j ) , $$ which, by the triangle inequality, implies  
 $$  | u_{i}^\top E U_{\perp} E u_{j}| \le | Eu_i \cdot Eu_j | + \sum_{k \le r} | (u_i^\top E u_k) ( u_k^\top E u_j )|. $$
\noindent By the definitions of $x$ and $y$ (see  \eqref{def: xyz}), 
$ | Eu_i \cdot Eu_j |  \le y,  \,\,\, | u_i^\top E u_j | \le x, $ for all $1 \le i \neq j \le r $. Therefore, 
 $$|u_i^\top EU_{\perp} E u_j| \leq  y + rx^2.$$
 Since $S_1$ contains only non-constant index sequences, every summand contains at least one adjacent pair $i_l \neq i_{l+1}$ for some $1 \leq l \leq k-1$. Thus, using the refined bound for one unequal adjacent pair and the trivial bound $\|E\|^2$ for the remaining $k-2$ factors, we obtain
\begin{equation} \label{S1product}
    \prod_{l=1}^{k-1} | u_{i_l}^\top E U_{\perp} E u_{i_{l+1}}| \leq (y+ rx^2) \|E\|^{2(k-2)}.
\end{equation}
Combining \eqref{S1integral}, \eqref{S1product}, and the fact that $\|u_{i_1}u_{i_k}^\top\| \leq 1$, we obtain
$$S_1 \leq \sum_{\substack{i_1,i_2,...,i_k \leq r \\ (i_1,...,i_k) \neq  (j,...,j) \\ \forall 1 \leq j \leq r}}  \frac{ 2^s \lambda_p}{ \lambda_p^{s/2} \delta_p^{ s/2 } } \times  (y+ rx^2) \|E\|^{2(k-2)}.$$
Since the number of summands in the definition of $S_1$ is at most $r^k$, we obtain 
$$S_1 \leq r^k \times \frac{ 2^s \lambda_p}{ \lambda_p^{s/2} \delta_p^{ s/2 } } \times  (y+ rx^2) \|E\|^{2(k-2)}.$$
 Replacing $k$ by $s/2+1$ and rearranging the RHS, we have 
 \begin{equation} \label{boundS1} S_1 \le 4 r \lambda_p \Big(  \frac{ ry +r^2 x^2}{ \lambda_p \delta_p }  \Big) \Big( \frac{2\sqrt{r}\|E\|}{\sqrt{\lambda_p \delta_p}} \Big)^{s-2} . 
 \end{equation} 

 \vskip2mm 

\noindent \textit{Bounding $S_2$.}  Since $\|u_j u_j^\top\| =1$ for all $1 \leq j \leq p,$ $S_2$ simplifies to 
$$ \sum_{j=1}^{p}  \norm{\frac{1}{2 \pi \textbf{i}}\int_{\Gamma} \frac{1}{z^{k-2} (z-\lambda_j)^{k}} dz}  \times  |  u_{j}^\top E U_{\perp} E u_{j}|^{k-1}.$$
As $ | u_{j}^\top E U_{\perp} E u_{j}|  \le \| E\| ^2$ for all $1 \leq j \leq p$, we have 
$$S_2 \leq  \sum_{j=1}^{p}  \norm{\frac{1}{2 \pi \textbf{i}}\int_{\Gamma} \frac{1}{z^{k-2} (z-\lambda_j)^{k}} dz}  \times  \|E\|^{2(k-1)}.$$
For each $1 \leq j \leq p$, to bound $ \norm{\frac{1}{2 \pi \textbf{i}}\int_{\Gamma} \frac{1}{z^{k-2} (z-\lambda_j)^{k}} dz}$, we apply Lemma \ref{key contour0} with the sets $X:=\{\lambda_j\}$ and $Y=\{0\}$. For this particular choice of $X$ and $Y$, $\lambda=\delta=\lambda_j$ and the total sum of exponents is $2k-2$.  Lemma \ref{key contour0} yields 
 \begin{equation*}
     \begin{split}
           &  \norm{ \frac{1}{2 \pi \textbf{i}}\int_{\Gamma} \frac{1}{z^{k-2} (z-\lambda_j)^{k}} dz }  \le \frac{2^{2k-2}}{\lambda_j^{2k-3}} \leq \frac{ 2^{2k-2}}{ \lambda_p^{2k-3} }.
     \end{split}
 \end{equation*}
 The last inequality follows from $\lambda_j \geq \lambda_p \geq 0.$ Therefore, 
 $$S_2 \leq  \sum_{j=1}^{p}  \norm{\frac{1}{2 \pi \textbf{i}}\int_{\Gamma} \frac{1}{z^{k-2} (z-\lambda_j)^{k}} dz}  \times  \|E\|^{2(k-1)} \leq p \frac{ 2^{2k-2}}{ \lambda_p^{2k-3} } \times \|E\|^{2(k-1)} = p \lambda_p \left( \frac{2\|E\|}{\lambda_p}\right)^{2k-2}.$$
 As $s=2k-2$, we finally have 
 \begin{equation} \label{boundS2} S_2 \le p \lambda_p \left( \frac{2 \| E\|}{ \lambda_p }  \right)^s.  
 \end{equation}

 Combining \eqref{PEQEPsplit}, \eqref{PEQEPSigmasplit}, \eqref{PEPEQS2split}, \eqref{boundS1}, and \eqref{boundS2}, we obtain 
 \begin{equation*}
     \begin{split}
   \bigg\|\frac{1}{2 \pi \textbf{i}}\int_{\Gamma}z PEQEPEQ \cdots PEQEP dz \bigg\| &\leq S_1+S_2    \\
   &  \leq   4r \lambda_p \Big(  \frac{ ry +r^2 x^2}{ \lambda_p \delta_p }  \Big) \Big( \frac{2\sqrt{r}\|E\|}{\sqrt{\lambda_p \delta_p}} \Big)^{s-2}+  p \lambda_p \left( \frac{2 \| E\|}{ \lambda_p }  \right)^s.
     \end{split}
 \end{equation*}
Dividing both sides by $\lambda_p$, we have 
\begin{equation} \label{badbound}
\frac{\Norm{\frac{1}{2 \pi \textbf{i}}\int_{\Gamma}M(\alpha; \beta) dz}}{\lambda_p} \leq    4r \Big(  \frac{ ry +r^2 x^2}{ \lambda_p \delta_p }  \Big) \Big( \frac{2\sqrt{r}\|E\|}{\sqrt{\lambda_p \delta_p}} \Big)^{s-2}+  p \left( \frac{2 \| E\|}{ \lambda_p }  \right)^s, 
\end{equation} for the special pair 
$(\alpha; \beta)$ when $T_1 + T_2 =0$. This proves the second part of Lemma \ref{lemma: Mab}.

We complete the whole proof. 
 %
 %

\section{An application: Low-rank approximation with missing and noisy entries} \label{sec: Apps}

\subsection{The problem and new results.} 
As outlined in the introduction, computing low-rank approximations of large data matrices is fundamentally important. In practice, however, data is often incomplete and noisy, prompting a central problem:

\begin{problem} \label{central} Compute a reasonably accurate low-rank approximation of 
the ground-truth matrix from imperfect data.  \end{problem}

Formally, let $A$ be the ground-truth matrix with entries 
$a_{ij}$, $1\le i \le m, 1 \le j \le n$. The noisy version is $A'$, with entries $a_{ij}' = a_{ij}+ \xi_{ij}$, where $\xi_{ij}$ are mean-zero independent (but not necessarily iid) random variables. 
We make the standard assumption that the set of observed entries is random. Technically, we assume that each entry of $A'$ is observed with probability $\rho$, independently.  For the unobserved entries, we write zero. The resulting matrix $B$, consisting of observed entries and zeroes elsewhere, is our input, from which we try to recover the best rank-$p$ approximation of $A$, for some given small parameter $p$.

Let $\tilde A := \rho^{-1} B $.  It is easy to see that $\E \tilde A = A$. Therefore, we can write 
$\tilde A = A +E$, where $E$ is a random matrix with independent entries having zero mean. For a given target rank $p$, we compute $\tilde A_p$, the best rank-$p$ approximation of $\tilde A$. Using our new results, we can prove that in many settings, $\tilde A_p $ is a good approximation for $A_p$. 
As an illustration, we analyze the following setting.

We assume that $m = \Theta (n)$ with $m \leq n$ and $r= \mathrm{rank}\, A=  O(1)$. Assume furthermore that 
$A$ has bounded entries,  $\| A \|_{\infty} = L = O(1)$, and 
$\|A\|_F^2 =\Theta (mn) =\Theta (n^2) $.
Notice that $\sum_{i=1} ^r \sigma_i(A)^2 =\|A\|_F^2= \Theta (n^2)$. We assume that the  $p$th largest singular value $\sigma_p =\Theta (n)$. Finally, assume  that the noise variables  $\xi_{ij}$
 are independent (but not necessarily iid), and each has mean 0 and bounded variance $\sigma_{ij}^2 \le \sigma^2=  O(1)$. 

These assumptions are imposed only to keep the presentation concise and avoid unnecessary technicalities.  Our argument extends to much more general settings, in line with the generalization of Theorems \ref{toy1}--\ref{toy3} (for which we used similar assumptions)  in Section \ref{sec: randomnoise}. The accuracy of this estimator is controlled by the following theorem.
\begin{theorem} \label{cor: fastcom} Assume the above setting.  If $\rho =\omega ( \frac{\log n}{n}) $, then with probability  $1 - o(1),$
$$ \textstyle \|\tilde A_p -A_p \| =  O \left( \sqrt{n/\rho}  \right)  + 
   \tilde O \left( \frac{n}{\rho \delta_p}+ \frac{ \sigma_p }{\sqrt{\rho} \delta_p} \right).$$
\end{theorem}

\subsection {Related results.}

One way to solve Problem \ref{central} is to first complete the matrix and then compute the approximation. 
The first step is the famous matrix completion problem; see Appendix \ref{App: motivating examples}. There are two cases:

{\bf Case 1.} If the matrix can be recovered fully, then we can compute the low-rank approximation exactly.  However, 
the assumptions needed for matrix completion,
particularly the incoherence assumption regarding the singular vectors, are strong and not always met in practice.

{\bf Case 2.} Without the incoherence assumption, one can still recover a matrix $M$ satisfying $ \| M  - A \|_F \le \Delta $. For a symmetric matrix $A$ of rank $r=O(1)$, the best results currently are \(\Delta= O\left( (\sigma_1/\sigma_r)^2 \cdot \sqrt{\frac{n}{\rho}}\right) \) via a gradient descent approach \cite{KMO2,chen2020nonconvex}, or $\Delta = O\left( (n\rho)^{1/4} \cdot \sqrt{\frac{n}{\rho}}\right)$ via a threshold approach \cite{Cha1}.

A natural way to proceed is to compute $M_p$ to approximate
$A_p$. It is not clear how to bound $\| M_p -A_p \|$. The standard argument using the Eckart--Young--Mirsky theorem would yield 
$$\|M_p -A_p\| \leq 2 \left(\sigma_{p+1}+ \|M-A\| \right) \leq 2 \left(\sigma_{p+1}+ \|M-A\|_F \right) \leq O \left(\sigma_{p+1} + \Delta \right).$$

Another relevant topic is denoising, for which there is an extensive literature; we refer to \cite{tong2025uniform} for an overview. The typical setting is that one assumes that $r$, the rank of $A$, is known, and tries to find a good rank $r$ approximation of $A$, using $B$ as the input. The standard assumption here is that $\sigma_r \ge 
C\| E \|$, for a sufficiently large constant $C$. Notice that in this case the gap $\delta_r= \sigma_r$, and many researchers have used the Davis-Kahan bound exploiting this large gap assumption; see, for instance, \cite{tong2025uniform, zhang2022heteroskedastic} and the references therein. The methods used in this topic do not seem to extend the best rank-$p$ 
approximation, when $p <r$, and the gap $\delta_p$ is not guaranteed to be larger than $\|E \|$.

Another closely related topic is accelerating computation using randomized sparsification;  see \cite{AchlioptasMcSherry2007, frieze2004fast, azar2001spectral, clarkson2017low}. Here, we have full access to the matrix $A$ and then sparsify it by keeping each entry with a small probability $\rho$ and setting the discarded entries to zero. The reason for doing so is that it is significantly faster to compute with the (rescaled) sparsified matrix $\tilde A$, as the majority of its entries are zero. The trade-off is in the accuracy of the output, precisely the problem we consider here.  However, when the goal is to compute a low-rank approximation, 
we are not aware of direct estimates for $\| \tilde A_p - A_p \|$. Researchers usually bound the weaker metrics
$$ \Big | \| A - A_p \| - \| \tilde A - \tilde A_p \|  \Big |; $$ see \cite{AchlioptasMcSherry2007, halko2011finding, frieze2004fast}. This quantity can be bounded by $\| E \| $ using the triangle inequality. This bound is often sharp, but does not always yield useful information for  $\| \tilde A_p -A_p \|$, as
the triangle inequality can be wasteful; see the discussions in the introduction and Subsection \ref{subsec: geoview}. 

\bibliographystyle{amsrefs}
\bibliography{RefO}

\appendix


 \section{ Proofs of the main theorems in Section \ref{sec: randomnoise} - Theorem \ref{theo: mainTh-1}, Theorem \ref{mainTh01}, and Theorem \ref{cor: AppTrivial}} \label{xy random}
 
 We first present the proof of Theorem \ref{mainTh01}. Theorem \ref{theo: mainTh-1} follows as a particular case. 

 Roughly speaking, Theorem \ref{mainTh01} is the random version of Theorem \ref{cor: rec}. It remains to estimate the skewness parameters $x,y$ given that $E$ is a $(K, \sigma^2)$-bounded random matrix. The idea is to use the moment method: we estimate the first and second moments of 
 $$\norm{u_i^\top E v_j} \,\,\text{for each}\,\,1 \leq i,j \leq r,$$
 and 
 $$\norm{E v_i \cdot E v_j}, \norm{E^\top u_i \cdot E^\top u_j}\,\,\text{for each}\,\, 1 \leq i < j \leq r,$$
 and then use the union bound over the choices of $(i,j)$. 
 
Therefore, given Theorem \ref{cor: rec}, Theorem \ref{mainTh01} is obtained via the following lemmas. The proofs of these lemmas and the moment computations can be found in \cite[Section 10 and Section 11]{DKTranVu1}.
\begin{lemma} \label{uEv} Let $E=(\xi_{kl})_{k,l=1}^{n}$ be a $m \times n$, $(K,\sigma^2)$-bounded random matrix. Then for any $1\leq i \leq j \leq r$ and every $t >0$,
$$\P( |u_i^\top E v_j| \geq t) \leq \exp \left( \frac{-t^2/2}{2\sigma^2+ Kt} \right).$$
\end{lemma}
 The next lemma provides an estimate for $y$. 
 \begin{lemma} \label{u_iEEu_j} Let $E=(\xi_{kl})_{k,l=1}^{n}$ be a $m \times n$, $(K,\sigma^2)$-bounded random matrix. Then, for any $t >0$, 
$$\max\left\lbrace \P (| E^\top u_1 \cdot E^\top u_2 - \mu_{12}  | \ge t),  \P (| E v_1 \cdot E v_2 -\mu'_{12}  | \ge t)\right\rbrace \le \frac{ 5(m+n) m_4 } {t^2} , $$ 
where 
$$\mu_{12}:= \E \left[ (E^\top u_1) \cdot (E^\top u_2) \right], \,\,\,\mu'_{12}:= \E \left[Ev_1 \cdot Ev_2  \right].$$
Using $\E \xi_{ij}^2 =\sigma_{ij}^2$ for $1 \leq i \leq m, 1 \leq j \leq n$, we can equivalently rewrite $\mu_{12}, \mu'_{12}$ as 
$$\mu_{12}= \sum_{j'=1}^{n} \sum_{j=1}^{m} \sigma_{j j'}^2 u_{2j} u_{1j} , \,\,\mu'_{12}= \sum_{i'=1}^{m} \sum_{i=1}^{n} \sigma_{i' i}^2 v_{2i} v_{1i}. $$

In particular, if the entries of $E$ are iid, then
$$\max\left\lbrace \P (| E^\top u_1 \cdot E^\top u_2 | \ge t),  \P (| E v_1 \cdot E v_2  | \ge t)\right\rbrace \le \frac{ 5(m+n) m_4 } {t^2} . $$ 

\end{lemma} 

\begin{proof}[Proof of Theorem \ref{cor: AppTrivial}]
To prove Theorem \ref{cor: AppTrivial}, we simply combine Theorem \ref{cor: rec} with Lemma \ref{uEv} and the trivial bound $y \leq \|E\|^2.$
\end{proof}

\section{Accurate low-rank approximation of the ground-truth matrix from imperfect data--Proof of Theorem \ref{cor: fastcom}}
We first recall the problem of computing a low-rank approximation of a ground-truth matrix from imperfect observations. Let $A=(a_{ij})\in\mathbb R^{m\times n}$ be the ground-truth matrix. Its noisy version is
\[
A'=(a'_{ij}),
\qquad
a'_{ij}=a_{ij}+\xi_{ij},
\]
where the $\xi_{ij}$ are independent, mean-zero random variables, not necessarily identically distributed.
Each entry of $A'$ is observed independently with probability $\rho$, producing the observed matrix $B$. Define
\[
\tilde A:=\rho^{-1}B.
\]
Then $\E\tilde A=A$, so we may write
\[
\tilde A=A+E,
\]
where $E$ is a random matrix with independent, mean-zero entries.

The goal of this section is to prove Theorem~\ref{cor: fastcom}, which provides a bound on
\[
\|\tilde A_p-A_p\|.
\]

In fact, Theorem~\ref{cor: fastcom} follows directly from Theorem~\ref{cor: AppTrivial} in this setting. Indeed, Theorem~\ref{cor: AppTrivial} gives
\[
\|\tilde{A}_p -A_p\| \leq 32 \bigg[p\|E\| + \frac{prx \cdot \sigma_p}{\delta_p} + \frac{r^2 \|E\|^2}{\delta_p} \bigg],
\]
where 
\[
x: =\max_{i,j \le r} | u_i^\top E v_j |.
\]
Thus, it remains to estimate $\|E\|$ and $x$. We use the following standard estimates from random matrix theory:
\begin{itemize}
\item   If  $\rho \gg \frac{\log n}{n}$, then $\|E\|=O \left( \sqrt {n/\rho}\right)$ with probability $1-o(1)$; see \cite{Vu0, bandeira2016sharp}.
\item  There is a constant $C$, depending on $L= \| A \|_{\infty} $, such that for  any fixed unit vectors $u,v$,
$$ \textstyle \P \left(|u^\top E v| \geq  \frac{ C }{\sqrt{\rho}} t \right) \leq \exp \left( -t^2 \right); $$ see \cite[Lemma 35]{OVK13}. It follows from the union bound that 
$$ \textstyle \P \left( x \geq  \frac{C} {\sqrt{\rho} } \log (m+n) \right) \leq r^2 (m+n)^{-1} =o(1), $$
where  $x: =\max_{i,j \le r} | u_i^\top E v_j | $ (see Definition \ref{def: xyz}). 
\end{itemize}
Substituting these estimates for $\|E\|$ and $x$ into Theorem~\ref{cor: AppTrivial} yields Theorem~\ref{cor: fastcom}.

We leave it to the reader to verify that this bound improves upon those obtained by classical methods, such as the Eckart--Young--Mirsky bound in Theorem~\ref{EY0} and the Davis--Kahan bound discussed in Remark~\ref{remark: daviskahan}.

Let us conclude this section with an illustrative numerical experiment (Figure \ref{fig: lowrank_approximate}):
\begin{figure}[ht]
    \centering
    \begin{subfigure}{0.45\textwidth}
        \includegraphics[width=\linewidth]{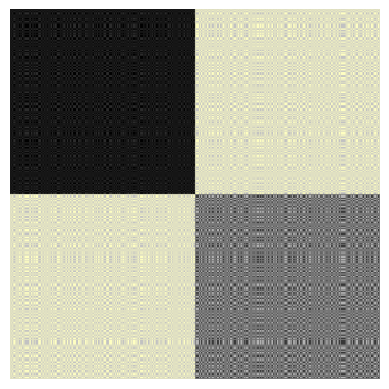}
        \caption{(a) Original Matrix $A$.}
    \end{subfigure}
    \hfill
    \begin{subfigure}{0.45\textwidth}
        \includegraphics[width=\linewidth]{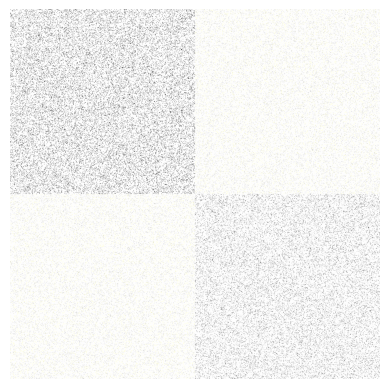}
        \caption{(b) Missing version $B$ with $\rho=0.1$.}
    \end{subfigure}

    \vspace{0.5em} 

    
    \begin{subfigure}{0.45\textwidth}
        \includegraphics[width=\linewidth]{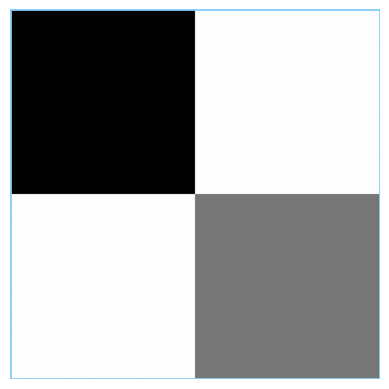} 
        \caption{(c) Rank-2 approximation of $A$.}
    \end{subfigure}
    \hfill
\begin{subfigure}{0.45\textwidth}
        \includegraphics[width=\linewidth]{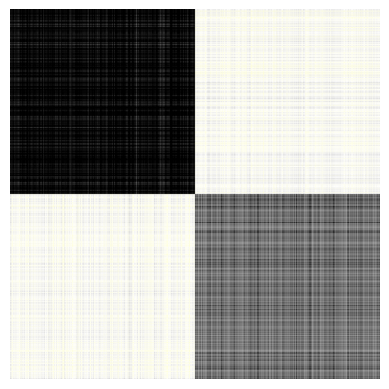}
        \caption{(d) Rank-2 approximation of $\rho^{-1}B$.}
    \end{subfigure}
    
    \caption{Visualizations of (a) the original matrix $A$ ($n=1000$, entries in $[-5,27]$), 
    (b) its missing version $B$ ($\rho=0.1$), 
    (c) the rank-2 approximation of $A$, and 
    (d)  the rank-2 approximation of $\rho^{-1}B$. Yellow indicates negative values, black positive, and white zero.}
    \label{fig: lowrank_approximate}
\end{figure}
 $A$ is a $1000 \times 1000$ matrix of rank $3$: 
\(
A = \sigma_1 u_1 u_1^\top + \sigma_2 u_2 u_2^\top + \sigma_3 u_3 u_3^\top,
\)
where  
\begin{itemize}
    \item $\sigma_1 = 22, \; \sigma_2 = 11, \; \sigma_3 = 5$;
    \item $u_1 = [\,\underbrace{1,\dots,1}_{500}, \underbrace{0,\dots,0}_{500}\,]^\top, \quad
            u_2 = [\,\underbrace{0,\dots,0}_{500}, \underbrace{1,\dots,1}_{500}\,]^\top$;
    \item $u_3$ has $500$ entries of value $-1$ and $500$ entries of value $1$, with random positions conditioned on being orthogonal to $u_1$ and $u_2$. 

    \item $\rho =.1$; namely, we observe ten percent of the entries. 

\end{itemize}

 \section{Motivating Examples} \label{App: motivating examples}
\subsection{Matrix recovery} \label{sec: matrix recovery}

A large matrix $A \in \mathbb{R}^{m\times n}$ is hidden, except for a few revealed entries in a set $\Omega \subset [m]\times [n]$.
We call $\Omega$ the set of \emph{observations} or \emph{samples}.
The matrix $A_\Omega$, defined by
\begin{equation}  \label{eq:sample-matrix}
    (A_\Omega)_{ij} = A_{ij}
    \ \text{ for } \
    (i, j)\in \Omega,
    \ \text{ and } 0
    \text{ otherwise},
\end{equation}

\noindent is called the \emph{observed} or \emph{sample} matrix.
The task is to recover $A$, given  $A_\Omega$. This is the \emph{matrix recovery} (or \emph{matrix completion}) problem, a central 
problem in data science that has received lots of attention in recent years, motivated by a number of real-life applications.
{
Examples include building recommendation systems, most notably the \textbf{Netflix challenge} \cite{netflix}, which has raised considerable attention from mainstream media; reconstructing a low-dimensional geometric surface based on partial distance measurements from a sparse network of sensors; repairing missing pixels in images; and system identification in control.
See the surveys \cite{li2019survey} and \cite{davenportSurvey2016} for more applications.}
In this discussion, we restrict ourselves to \emph{exact recovery}, where we want to recover all entries of $A$ exactly.

It is standard to assume that the set $\Omega$ is \emph{random}, and researchers have proposed two models: (a) $\Omega$ is sampled uniformly among subsets with the same size, or (b) that $\Omega$ has independently chosen entries, each with the same probability $p$, called the \emph{sampling density}, which can be known or hidden. It is often simple to replace the former model with the latter, using a simple conditioning trick. One samples the entries independently, and conditions on the event that the sample size equals a given number.

Let us focus on the independent entries model. Define 
$\tilde A= \frac{1}{p} A_{\Omega} $. Thanks to the normalization, this matrix  $\tilde A$ is a random matrix whose expectation is $A$. Thus, if we define $E := \tilde A - A$, then the (normalized) input matrix $\tilde A$ has the form
$$\tilde A = A +E,$$ where $E$ is a random matrix with independent entries having zero mean. The standard assumption on $A$ is that it has {\it low rank}  and delocalized eigenvectors (more precisely, $A$ satisfies the so-called incoherent condition) \cite{CT1, AV, recht2011simpler, hardt2014understanding, KMO1, AFW1, CR1} .

There are many approaches to attack the problem \cite{CT1, AV, recht2011simpler, hardt2014understanding, KMO1, AFW1, CR1,jain2013low,TranlinhVu}; 
we refer to \cite{TranlinhVu} for a survey. Many algorithms for this problem start with computing a low-rank approximation $B$ of $\tilde A$, and compare it with $A$
\cite{TranlinhVu, Cha1}. As $A$ itself has low rank, one hopes that the two matrices 
are close to each other, and then one can refine $B$ to get $A$ precisely. 
It is clear that estimation of the error term is critical here. 
In a recent paper, \cite{TranlinhVu}, L. Tran et al. applied the 
contour technique used in this paper to design a simple and fast algorithm that recovers $A$ precisely, even in the presence of noise, exploiting the fact that data in real life has a finite precision. 
(For instance, all entries of the Netflix matrix are half-integer.)

 \subsection{The spiked model in statistics }

 Let $X = [ \xi_1 \,\, \xi_2 \,\,...\,\, \xi_d]^\top$ be a random vector with  covariance matrix $A = (a_{ij})_{1 \leq i, j \leq d},$ where $a_{ij} = a_{ji} = \mathrm{Cov}(\xi_i, \xi_j).$  A fundamental problem in data science is estimating the unknown covariance matrix $A$ using the sample covariance matrix
 $$\tilde{A} := \frac{1}{n} \sum_{i=1}^{n} X_i X_i^\top,$$ where  $X_i$ are i.i.d samples of $X$.

 A primary task here is to approximate the leading eigenvectors and eigenspaces of $A$ from the input $\tilde{A}$. Among others, this is a central part of Principal Component Analysis (PCA), a statistical tool widely used in many scientific disciplines. Mathematically, the problem is essentially a perturbation problem for the pair $(A, \tilde{A})$. The matrix 
$$E = \tilde{A} - A$$ is a random matrix with mean $0$. Notice that the entries of $E$ are not independent, and their distributions vary with $X$.  On the other hand, it is still within the coverage of our general result. 

There is a huge literature on covariance matrices, and the spiked models have attracted significant attention in recent years; see, for instance, \cite{J1, K1, FW1, BYZ1} and the references therein. 
Informally speaking, in a spiked model,  the matrix  $A$ has a few important (large) eigenvalues (the "\textit{spikes}"),  while the remaining ones are small. This means that $A$ is approximately low rank. This case can be covered by an extension of our results here to such matrices; details will appear in a future paper \cite{DKTranVu2}.


 \subsection{Stochastic block model (hidden partition) } 
Statistical Block Model  (SBM) is a popular problem in statistics and combinatorics, where it is often referred to as the hidden partition problem. Here is the general setting.  Given a vertex set $V = \{1, \dots, n \}$, we partition $V$ into 
$k$ subsets $V_1, \dots, V_k$, where $k$ is small (typically a constant). Between a vertex  $u \in V_i$ and a vertex $v \in V_j$ ($1 \le i \le j \le k$), we draw a random edge with probability $p_{ij}$, independently. 
The goal is to recover the partition, under some conditions on the densities $p_{ij}$, given one instance of the resulting random graph. There is a substantial literature on this problem; see, for instance,  \cite{AKS1, BSp1, CK1, FR1, RCY1} and the references therein. 

Let $\tilde A$ be the adjacency matrix of the random graph and $A$ its expectation. We can write  $\tilde A= A +E$, where $E:= \tilde A -A $ is a (block) random matrix with zero mean and independent entries. 
By construction,  $A$ is the density matrix and has at most $k$ different rows. Therefore,  it has a rank at most $k$. 
Similar to the previous two problems,  the "noise" matrix $E$ does not represent white noise. We define it out of the context of the problem. 

The spectral method is a natural and popular approach to this problem; see \cite{AKS1, VuSpec1, AFW1,  RCY1}. Its success
is based on the fact that the leading eigenvectors of $A$ often reveal the partition. Since $A$ is unknown, one naturally computes the leading eigenvectors of the input $\tilde{A}$ and uses perturbation bounds to show that they are good approximations.
Next, one translates this approximation into an approximation of the partition; see  \cite{AKS1, AV, AFW1, VuSpec1}
for many arguments of this type.  
In many cases, one can further refine the approximate partition to obtain the exact partition; see  \cite{AKS1, AV, VuSpec1}. 
In this approach, a widely used assumption is that the gaps between the leading eigenvalues of $A$ are sufficiently large (typically of order $\sqrt n$ when the densities $p_{ij}$ are of order $\Theta (1)$). This is due to the use of the classical Davis-Kahan theorem. 

Our new perturbation results enable us to attack the problem under a much weaker gap assumption. Consequently, we can solve this problem in a wider range of parameters (which involve the densities $p_{ij}$ and the sizes  $| V_i|$ in the partition), see \cite[Section 5.3]{DKTranVu1} for more detailed examples,  including a generalization of the well-known hidden clique problem. In this generalization,  one can hide not one, but many cliques,  of about the same size.
Given a minimum gap (of order $\tilde O (1)$) between the sizes of different cliques, we are able to recover all of them using a simple spectral algorithm. It seems that previous arguments do not apply when the gap is this small. For some details, see \cite[Chapter 4 and Chapter 5]{Ver1book}.

 \subsection{Gaussian mixture and generalization }
 
Consider a large matrix $A \in \mathbb{R}^{m \times n}$, whose columns $A_1, A_2, \dots, A_n \in \mathbb{R}^m$ are one of $k$ \textit{centers}\ (vectors),  $\theta_1, \theta_2, \dots, \theta_k \in \mathbb{R}^m$,  based on a cluster assignment $z \in [k]^n$: 
$$A_i = \theta_{z_i} \,\,\text{for all} \,\,1 \leq i \leq n.$$

The cluster assignment $z$ is hidden, and our task is to recover it via a (noisy) observation $\tilde{A}= A+E$, in which $E$ consists of independent column-wise noises, $E_i$ for $1 \leq i \leq n$, and  
$$\tilde{A_i} = A_i + E_i \,\,\text{for all}\,\,1 \leq i \leq n.$$

In engineering terms, we receive $n$ noisy signals from $k$ different centers, and need to decide which signal belongs to 
which center. In many situations, the noise was assumed to be iid Gaussian, and the problem has been referred to as the Gaussian 
mixture problem. Other distributions have also been considered. As a matter of fact, it seems natural to assume that each center 
emits a different kind of noise. Thus, it makes sense to consider a generalization where the noise matrix $E$ does not have iid Gaussian columns, but random columns with different distributions.

Spectral clustering has been one of the most popular approaches to solve the mixture problem in high dimensions; see, for instance, \cite{ZZ1, AFW1}. The general idea is to first compute the best rank-$p$ approximation $\tilde{A}_p$ of $\tilde{A}$, for some small natural number $p$,  and then apply a classical clustering algorithm on the columns of $\tilde A_p$, viewed as points in the space.  Perturbation bounds on eigenvalues and eigenspaces are natural tools to analyze such an algorithm; see  \cite{ZZ1, AFW1}.

The standard assumption for spectral clustering to work is that $k$ is relatively small, order $O(\sqrt{\Delta)}$, for $\Delta:= \min_{i,j \in [k]} \Norm{\theta_i -\theta_j}$, and the eigenvalue gap $\delta_p:= \lambda_p -\lambda_{p+1}$ is large ($ \gg \Norm{E}$); see  \cite{ZZ1}. Equivalently, it means $A$ has low rank and a large eigenvalue gap. This is another important problem where the noisy matrix $\tilde A= A+E$ plays the role of the input, and perturbation bounds play an important role in the analysis of existing algorithms.

\section{Some classical perturbation bounds} \label{others}

This section recalls standard results referenced throughout our paper. 

\begin{theorem}[Eckart--Young--Mirsky bound~\cite{EY1}] \label{EY}
Let \( A, \tilde{A} \in \mathbb{R}^{m \times n} \), and let \( A_p \), \( \tilde{A}_p \) denote their respective best rank-$p$ approximations. Set \( E := \tilde{A} - A \). Then,
\[
\| \tilde{A}_p - A_p \| \le 2\left( \sigma_{p+1} + \| E \| \right),
\]
where \( \sigma_{p+1} \) is the $(p+1)$st singular value of \( A \).
\end{theorem}

\begin{theorem}[Weyl's inequality~\cite{We1}] \label{stanTheo: Weyl}
Let \( A, \tilde{A} \in \mathbb{R}^{m \times n} \), and define \( E:=\tilde{A} - A \). Then, for any \( 1 \le i \le n \),
\[
|\tilde{\lambda}_i - \lambda_i| \le \|E\| \quad \text{and} \quad |\tilde{\sigma}_i - \sigma_i| \le \|E\|,
\]
where \( \lambda_i, \tilde{\lambda}_i \) are the $i$th eigenvalues of \( A \) and \( \tilde{A} \), and \( \sigma_i, \tilde{\sigma}_i \) are the corresponding singular values.
\end{theorem}

\begin{theorem}[Davis-Kahan-Wedin sine theorem \cite{DKoriginal, W1}]
    \label{DKtheorem} Let \( A, \tilde{A} \in \mathbb{R}^{m \times n} \), and define \( E:=\tilde{A} - A \).  Denote the projections onto the subspace generated by the $p$-leading left singular vectors of $A, \tilde{A}$ by $\Pi_p, \tilde{\Pi}_p$ respectively. Then,  
\begin{equation}    \|\tilde{\Pi}_p - \Pi_p \| \leq \frac{ 2 \Norm{E}}{\sigma_p - \sigma_{p+1}}.
\end{equation}
\end{theorem}
\subsection{A perturbation bound of low-rank approximations via the Davis-Kahan theorem}\label{subsec: dktolowrank}
In this subsection, we derive an upper bound for $\|\tilde A_p-A_p\|$ using the Davis--Kahan--Wedin theorem.

Recall that $A_p=\Pi_p A$ and $\tilde A_p=\tilde\Pi_p\tilde A$. Since $\tilde A=A+E$, the triangle inequality gives 
\begin{equation*}
    \begin{split}
    \| \tilde A_p - A_p \| & = \| \tilde \Pi_p  \tilde A -\Pi_p A \|   \\
    & = \| \tilde{\Pi}_p A + \tilde{\Pi}_p E - \Pi_p A\| \\
    & = \|(\tilde{\Pi}_p -\Pi_p)A + (\tilde{\Pi}_p -\Pi_p) E + \Pi_p E\| \\
    & \le \| \tilde \Pi_p - \Pi_p  \|\cdot \| A \|  +  \| \tilde \Pi_p - \Pi_p \| \cdot \| E \| + 
\| \Pi_p E \| .
    \end{split}
\end{equation*}
Using Theorem \ref{DKtheorem}, we can bound the last formula by 
$$ 2 \frac{ \| E \| \sigma_1 }{\delta_p } + 2 \frac{ \| E\|^2} {\delta_p} 
+ \| \Pi_p E \|.$$
Since $\delta_p = \lambda_p -\lambda_{p+1} \le  \sigma_1$ and  $2 \|E\| \leq \delta_p$ (the requirement of Theorem \ref{DKtheorem}), we have $$\| \Pi_p E \| \leq \|E\| \leq \|E\| \times \frac{\sigma_1}{\delta_p}\,\, \text{and}\,\, 2 \frac{\|E\|^2}{\delta_p} \leq  \|E\| \frac{\sigma_1}{\delta_p}. $$
It follows that 
\begin{corollary} \label{corDK}
\begin{equation} \label{viaDK}
 \| \tilde A_p - A_p \| \leq 2 \frac{ \| E \| \sigma_1 }{\delta_p } + 2 \frac{ \| E\|^2} {\delta_p} 
+ \| \Pi_p E \| \leq 4 \frac{\|E\| \sigma_1}{\delta_p}.   
\end{equation} \end{corollary}
\end{document}